# AN OPERATOR APPROACH FOR MARKOV CHAIN WEAK APPROXIMATIONS WITH AN APPLICATION TO INFINITE ACTIVITY LÉVY DRIVEN SDES[1]

BY HIDEYUKI TANAKA AND ARTURO KOHATSU-HIGA[2]

*Mitsubishi UFJ Trust Investment Technology Institute Co., Ltd.
and Osaka University*

Weak approximations have been developed to calculate the expectation value of functionals of stochastic differential equations, and various numerical discretization schemes (Euler, Milshtein) have been studied by many authors. We present a general framework based on semigroup expansions for the construction of higher-order discretization schemes and analyze its rate of convergence. We also apply it to approximate general Lévy driven stochastic differential equations.

**1. Introduction.** Weak approximation problems play an important role in the numerical calculation of $E[f(X_t(x))]$ where $X_t(x)$ is the solution of the stochastic differential equation (SDE)

$$\begin{aligned}
X_t(x) = x &+ \int_0^t \tilde{V}_0(X_{s-}(x))\,ds + \int_0^t V(X_{s-}(x))\,dB_s \\
&+ \int_0^t h(X_{s-}(x))\,dY_s
\end{aligned} \quad (1.1)$$

with smooth coefficients $\tilde{V}_0 : \mathbf{R}^N \to \mathbf{R}^N, V = (V_1, \ldots, V_d), h : \mathbf{R}^N \to \mathbf{R}^N \otimes \mathbf{R}^d$ whose derivatives of any order ($\geq 1$) are bounded. Here $B_t$ is a $d$-dimensional standard Brownian motion and $Y_t$ is an $d$-dimensional Lévy process associ-

Received April 2008; revised September 2008.
[1]This article is partially based on the results of the first author which appear in his Master thesis.
[2]Supported in part by grants of the Japanese and Spanish governments.

ated with the Lévy triplet $(b, 0, \nu)$ satisfying the condition

$$\int_{\mathbf{R}_0^d} (1 \wedge |y|^p) \nu(dy) < \infty$$

for any $p \in \mathbf{N}$.

Our purpose is to find a discretization scheme $(X_t^{(n)}(x))_{t=0,T/n,\ldots,T}$ for given $T > 0$ such that

$$|E[f(X_T(x))] - E[f(X_T^{(n)}(x))]| \leq \frac{C(T, f, x)}{n^m}.$$

We denote briefly by $E[f(X_T(x))] - E[f(X_T^{(n)}(x))] = \mathcal{O}(1/n^m)$ the above situation, and say that $X_T^{(n)}$ is a $m$th-order discretization scheme for $X_t$ or that $X_T^{(n)}$ is an approximation scheme of order $m$. The Euler scheme is a first-order scheme, and has been studied by many researchers. Talay and Tubaro [19] show the first-order convergence of the Euler scheme and second-order convergence with the Romberg extrapolation for continuous diffusions. The fact that the convergence rate of the Euler scheme also holds for certain irregular functions $f$ under a Hörmander type condition has been proved by Bally and Talay [3] using Malliavin calculus. For the general Lévy driven case, the Euler–Maruyama scheme was first studied in Protter and Talay [18], see also Jacod and Protter [8] and Jacod et al. [7] (for smooth $f$). The Itô–Taylor (weak-Taylor) high-order scheme is a natural extension of the Euler scheme although is hard to simulate due to the use of multiple stochastic integrals. A discussion on the Itô–Taylor scheme with the Romberg extrapolation can be found in Kloeden and Platen [9].

In the continuous diffusion case, some new discretization schemes (also called Kusuoka type schemes) which are of order $m \geq 2$ without the Romberg extrapolation have been introduced by Kusuoka [11], Lyons and Victoir [13], Ninomiya and Victoir [16], Kusuoka, Ninomiya and Ninomiya [12], Ninomiya and Ninomiya [15] and Fujiwara [5, 6] ($m = 6$). The rate of convergence of these schemes is closely related to the stochastic Taylor expansion, or series expansion of exponential maps on a noncommutative algebra.

The actual simulation is carried out using (quasi) Monte Carlo methods. That is, one computes $\frac{1}{N} \sum_{i=1}^{N} f(X_T^{(n),i}(x))$ where $X_T^{(n),i}(x)$, $i = 1, \ldots, N$, denotes $N$ i.i.d. copies of $X_T^{(n)}(x)$. Therefore, using the law of large numbers, the final error $\frac{1}{N} \sum_{i=1}^{N} f(X_T^{(n),i}(x)) - E[f(X_T(x))]$ is of the order $O(\frac{1}{\sqrt{N}} + \frac{1}{n^m})$. Then the optimal asymptotic choice of $n$ is $O(n^m) = O(\sqrt{N})$.

The goal of the present article is twofold. First, we introduce a general framework to study weak approximation problems from the standpoint of operator (semigroup) expansions. That is given two processes that have equal semigroup expansions up to some order lead after composition to two processes that are closed in law. This goal is not new. In fact, using PDE techniques, Milshtein and Talay between others proved various weak approximation results. Although our proof is essentially the same it gives a new viewpoint that will help in defining new approximation schemes.

The next idea is to decompose the generator associated with (1.1) in (say) $d+2$ components where each component is associated with each component of the driving process (the whole Lévy process is considered as one component). Then we prove that if each of these components is approximated with an error of order $m+1$ then the composition gives an error of order $m$. In the particular case that each component can be characterized as the semigroup of a flow-type process then the composition leads to a composition-type approximation scheme.

Secondly, using the above strategy we provide approximations for solutions of (1.1). In particular, our approximations are valid for infinite activity Lévy processes $Y$. We prove that in fact, if one uses the Asmussen–Rosiński idea of approximating the jumps of size smaller than $\varepsilon$ with a Brownian motion and we only simulate one jump of size bigger than $\varepsilon$ per each time interval in the approximation is enough to provide a first-order approximation procedure. Furthermore we give the necessary estimate to determine $\varepsilon$ as a function of $n$. For this approximation, we found it better to decompose the generator in $d+4$ components.

This paper is organized as follows. In Section 2, we introduce the main example and the goal for the first part of this article in explicit mathematical terms. The general framework is introduced in Section 3. In Section 4 we give the results of convergence rates of numerical discretization schemes in the general framework. In Section 5, we give a general result that states how to recombine the approximations to coordinate processes in order to approximate the semigroup associated to (1.1). Finally, in Section 6 we approximate each coordinate process and in particular, we define approximation schemes for Lévy driven SDEs.

**2. Weak approximation problem.** In order to better understand the abstract formulation in Section 3, we introduce here our main example. Let $(Y_t)$ be a $d$-dimensional Lévy process characterized by Lévy–Khintchin formula

$$
(2.1) \quad \begin{aligned} &E[e^{i\langle \theta, Y_t \rangle}] \\ &= \exp t\bigg(i\langle \theta, b\rangle - \frac{\langle \theta, c\theta\rangle}{2} + \int_{\mathbf{R}_0^d}(e^{i\langle \theta, y\rangle} - 1 - i\langle \theta, \tau(y)\rangle)\nu(dy)\bigg),\end{aligned}
$$

where $b \in \mathbf{R}^d$, $c \in \mathbf{R}^d \otimes \mathbf{R}^d$ (symmetric, semi-positive definite) and $\nu$ is a Borel measure on $\mathbf{R}_0^d := \mathbf{R}^d \setminus \{0\}$ satisfying that

$$
(2.2) \quad \int_{\mathbf{R}_0^d}(1 \wedge |y|^p)\nu(dy) < \infty.
$$

This measure $\nu$ is called the Lévy measure. It is well known that (2.2) implies that $Y_t \in \bigcap_{p \geq 1} L^p$ for all $t$. We also recall that $\tau$ is a truncation function [e.g., $\tau(y) = y 1_{\{|y| \leq 1\}}$, the constant $b$ and $\tau$ depend on each other]. The triplet $(b, c, \nu)$ is called the Lévy triplet.

The Lévy driven stochastic differential equation is given by

$$
\begin{aligned}
X_t(x) = x &+ \int_0^t \tilde{V}_0(X_{s-}(x)) \, ds + \int_0^t V(X_{s-}(x)) \, dB_s \\
&+ \int_0^t h(X_{s-}(x)) \, dY_s
\end{aligned}
\tag{2.3}
$$

with smooth coefficients $\tilde{V}_0 : \mathbf{R}^N \to \mathbf{R}^N, V = (V_1, \ldots, V_d), h : \mathbf{R}^N \to \mathbf{R}^N \otimes \mathbf{R}^d$ whose derivatives of any order ($\geq 1$) are bounded. Here $B_t$ and $Y_t$ are independent $d$-dimensional standard Brownian motion and $Y_t$ is a $d$-dimensional Lévy process associated with the Lévy triplet $(b, 0, \nu)$ satisfying condition (2.2). Using general semimartingale theory (see [17]) we have that the above equation has a unique solution. We define $V_0 := \tilde{V}_0 - \frac{1}{2} \sum_{i=1}^d \sum_{j=1}^N \frac{\partial V_i}{\partial x_j} V_i^{(j)}$. Then (2.3) can be rewritten in the following Stratonovich form:

$$
X_t(x) = x + \sum_{i=0}^d \int_0^t V_i(X_{s-}(x)) \circ dB_s^i + \int_0^t h(X_{s-}(x)) \, dY_s,
$$

where $B_t^0 = t$.

Before introducing the general framework of approximation, let us explain in mathematical terms the goal in this article. Our main example corresponds to the approximation of the semigroup $P_t$ defined as the semigroup associated to the Markov process $X_t$:

$$
P_t f(x) = E[f(X_t(x))],
$$

where $f : \mathbf{R}^N \to \mathbf{R}$ is a continuous function with polynomial growth at infinity.

Let $Q_t \equiv Q_t^n$ be an operator such that the semigroup property is satisfied in $\{kT/n; k = 0, \ldots, n\}$. Assume that $Q_t$ approximates $P_t$ in the sense that it satisfies the local error estimate $(P_t - Q_t)f(x) = \mathcal{O}(t^{m+1})$. Then using the semigroup property of both $P_t$ and $(Q_{kT/n})$, we notice that

$$
P_T f(x) - (Q_{T/n})^n f(x) = \sum_{k=0}^{n-1} (Q_{T/n})^k (P_{T/n} - Q_{T/n}) P_{T-(k+1)/nT} f(x).
$$

Therefore if we have good norm estimates of $(Q_{T/n})^k$ and $P_{T-(k+1)/nT}$ in a sense to be defined later (in particular the norm estimates have to be independent of $n$) then we can expect that $(Q_{T/n})^n$ is an approximation of order $m$ to $P_T$. Finally in order to be able to perform Monte Carlo simulations we assume that $Q$ has a *stochastic representation*. That is, there exists a stochastic process $M = M_t(x)$ starting at $x$ such that $Q_t f(x) = E[f(M_t(x))]$. Then clearly, we have the following representation.

$$Q_T f(x) = (Q_{T/n})^n f(x) = E[f(M_{T/n}^1 \circ \cdots \circ M_{T/n}^n(x))],$$

where $M_{T/n}^i$ are independent copies of $M_{T/n}$ and $\circ$ is defined as $(M_t^i \circ M_t^j)(x) := M_t^i(M_t^j(x))$.

The above ideas are well known and have been already used to achieve proofs of weak convergence (for historical references, see [9]). Nevertheless, it seems this is the first time it appears in this general framework. For example, if we take $M_t(x) := x + \tilde{V}_0(x)t + V(x)B_t + h(x)Y_t$ for $d=1$, one obtains the Euler–Maruyama scheme.

Next to further simplify the procedure to obtain approximations we write the operator $P_t$ as a composition of $d+2$ operators as follows. First define the following stochastic processes $X_{i,t}(x)$, $i=0,\ldots,d+1$, usually called coordinate processes, which are the unique solutions of

$$X_{0,t}(x) = x + \int_0^t V_0(X_{0,s}(x))\,ds,$$

$$X_{i,t}(x) = x + \int_0^t V_i(X_{i,s}(x)) \circ dB_s^i, \qquad 1 \leq i \leq d,$$

$$X_{d+1,t}(x) = x + \int_0^t h(X_{d+1,s-}(x))\,dY_s.$$

Then we define

(2.4) $$Q_{i,t} f(x) := E[f(X_{i,t}(x))]$$

for continuous function $f : \mathbf{R}^N \to \mathbf{R}$ with polynomial growth at infinity.

For notational convenience we identify a smooth function $V : \mathbf{R}^N \to \mathbf{R}^N$ with a smooth vector field $\sum_{i=1}^N V^{(i)} \frac{\partial}{\partial x_i}$ on $\mathbf{R}^N$. Let us define integro-differential operators $L_i$ acting on $C^2$ by

$$L_0 f(x) := (V_0 f)(x), \qquad L_i f(x) := \tfrac{1}{2}(V_i^2 f)(x), \qquad 1 \leq i \leq d,$$

(2.5) $$L_{d+1} f(x) := \nabla f(x) h(x) b$$
$$+ \int (f(x + h(x)y) - f(x) - \nabla f(x) h(x) \tau(y)) \nu(dy).$$

It is well known that $L := \sum_{i=0}^{d+1} L_i$ is the generator of $X$ and similarly $L_i$ is the generator of $X_{i,t}$. Also $e^{tL} := P_t$ and $e^{tL_i} := Q_{i,t}$, respectively, where we consider these expressions as exponential maps on a noncommutative algebra. One notices that these operators have the form

(2.6) $$e^{tL} = \sum_{k=0}^m \frac{t^k}{k!} L^k + \mathcal{O}(t^{m+1}),$$

$$(2.7) \quad e^{tL_i} = \sum_{k=0}^{m} \frac{t^k}{k!} L_i^k + \mathcal{O}(t^{m+1}).$$

To approximate $e^{tL}$, we would like to find some combination of operators satisfying

$$(2.8) \quad e^{tL} - \sum_{j=1}^{k} \xi_j e^{t_{1,j} A_{1,j}} \cdots e^{t_{\ell_j,j} A_{\ell_j,j}} = \mathcal{O}(t^{m+1})$$

with some $t_{i,j} > 0$, $A_{i,j} \in \{L_0, L_1, \ldots, L_{d+1}\}$ and weights $\{\xi_j\} \subset [0,1]$ with $\sum_{j=1}^{k} \xi_j = 1$. This will correspond to an $m$th-order discretization scheme.

To find such schemes, one can perform formal Taylor expansions for $e^{tA}$ in each of the terms in (2.8). We remark that the result (2.8) will follow directly from (2.6) and (2.7) independent of the specific form of the decomposition $L := \sum_{i=0}^{d+1} L_i$. This algebraic calculation has lead to the introduction of the following approximation schemes:

*Ninomiya–Victoir* (a):

$$(2.9) \quad \tfrac{1}{2} e^{t/2L_0} e^{tL_1} \cdots e^{tL_{d+1}} e^{t/2L_0} + \tfrac{1}{2} e^{t/2L_0} e^{tL_{d+1}} \cdots e^{tL_1} e^{t/2L_0}.$$

*Ninomiya–Victoir* (b):

$$(2.10) \quad \tfrac{1}{2} e^{tL_0} e^{tL_1} \cdots e^{tL_{d+1}} + \tfrac{1}{2} e^{tL_{d+1}} \cdots e^{tL_1} e^{tL_0}.$$

*Splitting method*:

$$(2.11) \quad e^{t/2L_0} \cdots e^{t/2L_d} e^{tL_{d+1}} e^{t/2L_d} \cdots e^{t/2L_0}.$$

The semigroups generated by these operators have a probabilistic representation. For example, Ninomiya–Victoir (a) corresponds to

$$1_{U<1/2} X_{0,t/2} \circ X_{d+1,t} \cdots X_{1,t} \circ X_{0,t/2}(x)$$
$$+ 1_{1/2 \leq U} X_{0,t/2} \circ X_{1,t} \cdots X_{d+1,t} \circ X_{0,t/2}(x),$$

where $U$ is a uniform random variable taking values in $[0,1]$, independent of $X_{i,t}$. However, since a closed-form solution $X_{i,t}$ is not always available, one has to replace $X_{i,t}$ with other approximations of order $m+1$ so that the final approximation result remains unchanged. Nevertheless the fact that there is only one driving process simplifies this task. This problem will be discussed in Section 5.

### 3. Preliminaries.

3.1. *Notation and assumptions.* In this section, we consider a general framework for weak approximations following the arguments in Section 2, without using the specific form of the operator. We first define the following functional spaces.

- $C_p^m \equiv C_p^m(\mathbf{R}^N)$: the set of $C^m$ functions $f:\mathbf{R}^N \to \mathbf{R}$ such that for each multi-index $\alpha$ with $0 \leq |\alpha| \leq m$, $|\partial_x^\alpha f(x)| \leq C(\alpha)(1+|x|^p)$ for some positive constant $C(\alpha)$.

  We also let $C_p \equiv C_p^0$. Let us define a norm on $C_p^m$ by
  $$\|f\|_{C_p^m} := \inf\{C \geq 0 : |\partial_x^\alpha f(x)| \leq C(1+|x|^p), 0 \leq |\alpha| \leq m, x \in \mathbf{R}^N\},$$
  where we denote $|\alpha| := \alpha_1 + \cdots + \alpha_N$ for $\alpha = (\alpha_1, \ldots, \alpha_N) \in \mathbf{Z}_+^N$.

- $C_p^{1,m}([0,T] \times \mathbf{R}^N)$: the set of functions $f:[0,T] \times \mathbf{R}^N \to \mathbf{R}$ such that $s \mapsto f(s,x)$ is continuous differentiable for all $x \in \mathbf{R}^N$ and satisfies that $f(s,\cdot)$, $\partial_s f(s,\cdot) \in C_p^m$ with $\sup_{s \in [0,T]}(\|f(s,\cdot)\|_{C_p^m} + \|\partial_s f(s,\cdot)\|_{C_p^m}) < \infty$.

From now on, we denote by $Q_t : \bigcup_{p \geq 0} C_p(\mathbf{R}^N) \to \bigcup_{p \geq 0} C_p(\mathbf{R}^N)$ a linear operator for $0 \leq t \leq T$.

ASSUMPTION $(\mathcal{M}_0)$. If $f \in C_p$ with $p \geq 2$, then $Q_t f \in C_p$ and
$$\sup_{t \in [0,T]} \|Q_t f\|_{C_p} \leq K \|f\|_{C_p}$$
for some constant $K > 0$ independent of $n$. Futhermore, we assume $0 \leq Q_t f(x) \leq Q_t g(x)$ whenever $0 \leq f \leq g$.

We now introduce two assumptions which are highly related to the convergence rate of approximation schemes.

ASSUMPTION $(\mathcal{M})$. $Q_t$ satisfies $(\mathcal{M}_0)$, and for each $f_p(x) := |x|^{2p}$ ($p \in \mathbf{N}$),

(3.1) $$Q_t f_p(x) \leq (1+Kt)f_p(x) + K't$$

for some constant $K = K(T,p)$, $K' = K'(T,p) > 0$.

For $m \in \mathbf{N}$, $\delta_m : [0,T] \to \mathbf{R}_+$ denotes a decreasing function which satisfies
$$\limsup_{t \to 0+} \frac{\delta_m(t)}{t^{m-1}} = 0.$$
Usually, we have $\delta_m(t) = t^m$.

ASSUMPTION $\mathcal{R}(m, \delta_m)$. For each $p \geq 2$, there exists a constant $q = q(m, p) \geq p$ and linear operators $e_k : C_p^{2k} \to C_{p+2k}$ ($k = 0, 1, \ldots, m$) such that

(A): For every $f \in C_p^{2(m'+1)}$ with $1 \leq m' \leq m$, the operator $Q_t$ satisfies

$$(3.2) \qquad Q_t f(x) = \sum_{k=0}^{m'} (e_k f)(x) t^k + (\mathrm{Err}_t^{(m')} f)(x), \qquad t \in [0, T],$$

where $\mathrm{Err}_t^{(m')} f \in C_q$, and satisfies the following condition:

(B): If $f \in C_p^{m''}$ with $m'' \geq 2k$, then $e_k f \in C_{p+2k}^{m''-2k}$ and there exists a constant $K = K(T, m) > 0$ such that

$$(3.3) \qquad \|e_k f\|_{C_{p+2k}^{m''-2k}} \leq K \|f\|_{C_p^{m''}}, \qquad k = 0, 1, \ldots, m.$$

Furthermore if $f \in C_p^{m''}$ with $m'' \geq 2m' + 2$,

$$\|\mathrm{Err}_t^{(m')} f\|_{C_q} \leq \begin{cases} K t^{m'+1} \|f\|_{C_p^{m''}}, & \text{if } m' < m, \\ K t \delta_m(t) \|f\|_{C_p^{m''}}, & \text{if } m' = m \end{cases}$$

for all $0 \leq t \leq T$.

(C): For every $0 \leq k \leq m$ and $j \geq 2k + 2$, if $f \in C_p^{1,j}([0, T] \times \mathbf{R}^N)$, then $e_k f \in C_{p+2k}^{1, j-2k}([0, T] \times \mathbf{R}^N)$.

In order to compare the finite power expansions of different operators, we introduce the following notation:

$$J_{\leq m}(Q_t) := \sum_{k=0}^{m} t^k e_k,$$

$$J_m(Q) := e_m.$$

$J_{\leq m}(Q_t)$ is a linear operator, which is related to the series expansion of $t \mapsto e^{tL_i}$ (cf. Proposition A.6). The following lemma comprises some basic properties related to the above definition. The proof is straightforward.

LEMMA 3.1. *The following properties are satisfied:*

$$\mathcal{R}(m+1, \delta_{m+1}) \Rightarrow \mathcal{R}(m, t^m),$$

$$\mathcal{R}(m, \delta_m) \Rightarrow \mathcal{R}(m, \tilde{\delta}_m),$$

whenever $\delta_m(t) \leq K \tilde{\delta}_m(t)$ and $\limsup_{t \to 0+} \tilde{\delta}_m(t)/t^{m-1} = 0$.

(i) *Let $\{\xi_i\}_{1 \leq i \leq \ell}$ be deterministic positive constants with $\sum_i \xi_i = 1$, and assume (M) for $Q_t^{(i)}$ ($i = 1, \ldots, \ell$). Then $\sum_{i=1}^{\ell} \xi_i Q_t^{(i)}$ also satisfies (M).*

(ii) *Let $\{\xi_i\}_{1 \leq i \leq \ell} \subset \mathbf{R}$ and assume $\mathcal{R}(m, \delta_m)$ for $Q_t^{(i)}$ ($i = 1, \ldots, \ell$). Then $\sum_{i=1}^{\ell} \xi_i Q_t^{(i)}$ also satisfies $\mathcal{R}(m, \delta_m)$.*

**4. Weak rate of convergence.** In this section, we prove the rate of convergence for the approximating operator $Q$ under the Assumptions $(\mathcal{M})$, $\mathcal{R}(m, \delta_m)$. Throughout this section, we assume the following assumption:

ASSUMPTION $(\mathcal{M}_P)$. For all $f \in C_p^m$ then $P.f \in C_{p+2}^{1,m-2}$ and furthermore the following two properties are satisfied for some positive constant $C$:
1. $\sup_{t \in [0,T]} \|P_t f\|_{C_p^m} \leq C\|f\|_{C_p^m}$,
2. $\|(P_t - P_s)f\|_{C_p^m} \leq C|t-s|\|f\|_{C_p^m}$.

THEOREM 4.1. *Assume $(\mathcal{M})$ and $\mathcal{R}(m, \delta_m)$ for $P_t$ and $Q_t$ with $J_{\leq m}(P_t - Q_t) = 0$. Then for any $f \in C_p^{2(m+1)}$, there exists a constant $K = K(T, x) > 0$ such that*

$$(4.1) \qquad |P_T f(x) - (Q_{T/n})^n f(x)| \leq K \delta_m\left(\frac{T}{n}\right) \|f\|_{C_p^{2(m+1)}}.$$

For the proof, we need the following lemma.

LEMMA 4.1. *Under Assumption $(\mathcal{M})$, the operators $P_t$ and $Q_t$ satisfy*

$$\sup_n \max_{0 \leq k \leq n} ((P_{T/n})^k + (Q_{T/n})^k) f(x) < \infty$$

*for any positive function $f \in C_p$ with $p \geq 0$.*

PROOF. Without loss of generality we do the proof for $Q$. Let $f_p(x) = |x|^{2p}$ for $p \in \mathbf{N}$. By the Assumption $(\mathcal{M})$, we have

$$(Q_{T/n})^k f_p(x) = (Q_{T/n})^{k-1}(Q_{T/n} f_p)(x)$$
$$\leq \left(1 + \frac{C}{n}\right)(Q_{T/n})^{k-1} f_p(x) + \frac{C'}{n}$$

with some constant $C, C'$ independent of $t, x, k, n$. Since $(1 + \frac{C}{n})^k \leq e^C$, one proves by induction that

$$\sup_n \max_{0 \leq k \leq n} (Q_{T/n})^k f_p(x) \leq e^C C'(1 + |x|^{2p}).$$

This completes the proof. $\square$

PROOF OF THEOREM 4.1. Let $f \in C_p^{2(m+1)}$. Using the semigroup property and Assumption $\mathcal{R}(m, \delta_m)$, we have

$$P_T f(x) - (Q_{T/n})^n f(x) = \sum_{k=0}^{n-1} (Q_{T/n})^k (P_{T/n} - Q_{T/n}) P_{T-(k+1)/nT} f(x)$$
$$= \sum_{k=0}^{n-1} (Q_{T/n})^k (\text{Err}_{T/n}^{(m)} P_{T-(k+1)/nT} f)(x),$$

where $\mathrm{Err}_t^{(m)}$ is the error term of $(P - Q)$ defined in (3.2).

We obtain from Assumptions $\mathcal{R}(m, \delta_m)$ and $(\mathcal{M}_P)$

$$|(\mathrm{Err}_{T/n}^{(m)} P_{T-(k+1)/nT} f)(x)| \leq K_1 \frac{T}{n} \delta_m\left(\frac{T}{n}\right)(1 + |x|^q) \|P_{T-(k+1)/nT} f\|_{C_p^{2(m+1)}}$$

$$\leq \frac{K_2 T}{n} \delta_m\left(\frac{T}{n}\right)(1 + |x|^q) \|f\|_{C_p^{2(m+1)}}$$

and hence Lemma 4.1 leads to

$$|(Q_{T/n})^k (\mathrm{Err}_{T/n}^{(m)} P_{T-(k+1)/nT} f)(x)|$$

$$\leq \frac{K_2 T}{n} \delta_m\left(\frac{T}{n}\right) \|f\|_{C_p^{2(m+1)}} (Q_{T/n})^k (1 + |x|^q)$$

$$\leq \frac{K}{n} \delta_m\left(\frac{T}{n}\right) \|f\|_{C_p^{2(m+1)}}$$

for some constant $K = K(T, x)$. This completes the proof. $\square$

The following theorem is an extension of Theorem 4.1, and is analogous to Talay and Tubaro [19], Theorem 1.

THEOREM 4.2. *Assume* $(\mathcal{M})$ *and* $\mathcal{R}(m+1, \delta_{m+1})$ *for* $Q_t$ *with* $J_{\leq m}(P_t - Q_t) = 0$. *Then for each* $f \in C_p^{2(m+3)}$, *we have*

$$(4.2) \quad P_T f(x) - (Q_{T/n})^n f(x) = \frac{K}{n^m} + \mathcal{O}\left(\left(\frac{T}{n}\right)^{m+1} \vee \delta_{m+1}\left(\frac{T}{n}\right)\right),$$

*where* $K = T^m \int_0^T P_s J_{m+1}(P - Q) P_{T-s} f(x)\,ds$.

PROOF. We start by noting that as in the proof of Theorem 4.1,

$$(P_{T/n} - Q_{T/n}) P_{T-s} f(x)$$

$$= \left(\frac{T}{n}\right)^{m+1} J_{m+1}(P - Q) P_{T-s} f(x) + (\mathrm{Err}_{T/n}^{(m+1)} P_{T-s} f)(x)$$

and therefore,

$$P_T f(x) - (Q_{T/n})^n f(x)$$

$$= \left(\frac{T}{n}\right)^{m+1} \sum_{k=0}^{n-1} (Q_{T/n})^k J_{m+1}(P - Q) P_{T-(k+1)/nT} f(x)$$

$$+ \mathcal{O}\left(\delta_{m+1}\left(\frac{T}{n}\right)\right).$$

Now applying the proof of Theorem 4.1 (for $m = 1$) to $J_{m+1}(P - Q) \times P_{T-(k+1)/nT}f \in C^4_{p+2(m+1)}$, we obtain

$$|((Q_{T/n})^k - P_{kT/n})J_{m+1}(P - Q)P_{T-(k+1)/nT}f(x)|$$
$$\leq \frac{C_1(T,x)}{n}\|J_{m+1}(P - Q)P_{T-(k+1)/nT}f\|_{C^4_{p+2(m+1)}}$$
$$\leq \frac{C_2(T,x)}{n}\|f\|_{C^{2(m+3)}_p}.$$

Next, we have by hypothesis $(\mathcal{M}_P)$,

$$|P_{kT/n}J_{m+1}(P - Q)P_{T-(k+1)/nT}f(x)$$
$$- P_{(k+1)/nT}J_{m+1}(P - Q)P_{T-(k+1)/nT}f(x)|$$
$$= |(I - P_{T/n})P_{kT/n}J_{m+1}(P - Q)P_{T-(k+1)/nT}f(x)|$$
$$\leq \frac{C_3(T,x)}{n}\|P_{kT/n}J_{m+1}(P - Q)P_{T-(k+1)/nT}f\|_{C^4_{p+2(m+1)}}$$
$$\leq \frac{C_4(T,x)}{n}\|f\|_{C^{2(m+3)}_p}.$$

Using Lemmas A.1 and A.2 in the Appendix and $J_{m+1}(P-Q)P_{T-s}f(x) \in C^{1,2}_{p+2(m+2)}$, we have

$$\left|\frac{T}{n}\sum_{k=0}^{n-1} P_{(k+1)/nT}J_{m+1}(P - Q)P_{T-(k+1)/nT}f(x)\right.$$
$$\left. - \int_0^T P_s J_{m+1}(P - Q)P_{T-s}f(x)\,ds\right|$$
$$\leq \frac{C(T,f,x)}{n}.$$

As a result, taking $K = T^m \int_0^T P_s J_{m+1}(P-Q)P_{T-s}f(x)\,ds$, we conclude that

$$P_T f(x) - (Q_{T/n})^n f(x) = \frac{K}{n^m} + \mathcal{O}\left(\left(\frac{T}{n}\right)^{m+1} \vee \delta_{m+1}\left(\frac{T}{n}\right)\right).$$

This concludes the proof. □

**5. Algebraic approximations of semigroup operators using coordinate operators.** Throughout this section, we assume that $P_t$, $t \in [0, T]$, is a semigroup that satisfies $(\mathcal{M})$, $(\mathcal{M}_P)$ and $\mathcal{R}(m, \delta_m)$. Furthermore we suppose that

$$J_{\leq m}(P_t) = I + \sum_{j=1}^{m} \frac{t^j}{j!} e_j$$

with $e_j = (\sum_{i=0}^{d+1} L_i)^j$ satisfying the properties stated in $\mathcal{R}(m, \delta_m)$. Similarly, we assume that $Q_{i,t} : \bigcup_{p \geq 0} C_p(\mathbf{R}^N) \to \bigcup_{p \geq 0} C_p(\mathbf{R}^N)$, $i = 0, \ldots, d+1$, be a sequence of operators such that they satisfy $(\mathcal{M})$, $(\mathcal{M}_P)$ and $\mathcal{R}(m, \delta_m)$ with

$$J_{\leq m}(Q_{i,t}) = I + \sum_{j=1}^{m} \frac{t^j}{j!} L_i^j.$$

$\prod_{i=1}^{\ell} a_i := a_1 a_2 \cdots a_\ell$ denotes a noncommutative product.

THEOREM 5.1. *Assume $m = 2$. That is, $(\mathcal{M})$ and $\mathcal{R}(2, \delta_2)$ are satisfied for $Q_{i,t}$ ($i = 0, 1, \ldots, d+1$). Then all the following operators satisfy $(\mathcal{M})$ and $\mathcal{R}(2, \delta_2)$:*

*N–V (a):* $Q_t^{(a)} = \frac{1}{2} Q_{0,t/2} \prod_{i=1}^{d+1} Q_{i,t} Q_{0,t/2} + \frac{1}{2} Q_{0,t/2} \prod_{i=1}^{d+1} Q_{d+2-i,t} Q_{0,t/2}$.

*N–V (b):* $Q_t^{(b)} = \frac{1}{2} \prod_{i=0}^{d+1} Q_{i,t} + \frac{1}{2} \prod_{i=0}^{d+1} Q_{d+1-i,t}$.

*Splitting:* $Q_t^{(sp)} = Q_{0,t/2} \cdots Q_{d,t/2} Q_{d+1,t} Q_{d,t/2} \cdots Q_{0,t/2}$.

*Moreover, we have $J_{\leq 2}(Q_t^{(a)}) = J_{\leq 2}(Q_t^{(b)}) = J_{\leq 2}(Q_t^{(sp)}) = \sum_{k=0}^{2} \frac{t^k}{k!} L^k$. In particular, the above schemes define a second-order approximation scheme.*

The proof of Theorem 5.1 is an application of Theorem 4.1. The conditions follow from the next lemma, together with an algebraic calculation as pointed out at the end of Section 2.

This theorem can also be stated for third-order approximation schemes.

LEMMA 5.1. *Let $Q_t^1$ and $Q_t^2 : \bigcup_{p \geq 0} C_p(\mathbf{R}^N) \to \bigcup_{p \geq 0} C_p(\mathbf{R}^N)$ be two linear operators and let $Q_t^1 Q_t^2$ be the composite operator. Then:*

(i) *If $(\mathcal{M})$ holds for $Q_t^1$, $Q_t^2$, then it also holds for $Q_t^1 Q_t^2$.*
(ii) *If $\mathcal{R}(m, \delta_m)$ holds for $Q_t^1$, $Q_t^2$, then it also holds for $Q_t^1 Q_t^2$.*

PROOF. (i) is obvious. We now prove (ii). Let $m' \leq m$. We have by hypothesis that

$$Q_t^1 f(x) = \sum_{k=0}^{m'} (J_k Q_t^1 f)(x) t^k + (\text{Err}_t^{(m',1)} f)(x),$$

$$Q_t^2 f(x) = \sum_{k=0}^{m'} (J_k Q_t^2 f)(x) t^k + (\text{Err}_t^{(m',2)} f)(x)$$

for $f \in C_p^{2(m'+1)}$, $p \geq 2$. Furthermore there exists $q = q(m, p) > 0$ such that $\text{Err}_t^{(m',1)} f$, $\text{Err}_t^{(m',2)} f \in C_q$. Now we prove (A)–(C) in the definition of $\mathcal{R}(m, \delta_m)$.

(A): Note that for $f \in C_p^{2(m'+1)}(\mathbf{R}^N)$,

$$Q_t^1 Q_t^2 f(x) = Q_t^1 \bigg( \sum_{k=0}^{m'} (J_k Q_t^2 f)(x) t^k + (\mathrm{Err}_t^{(m',2)} f)(x) \bigg).$$

Since $J_k Q_t^2 f \in C_{p+2k}^{2(m'+1)-2k}$, $Q_t^1(J_k Q_t^2 f)$ can be written as

$$(Q_t^1(J_k Q_t^2 f))(x) = \sum_{\ell=0}^{m'-k} (J_\ell Q_t^1 (J_k Q_t^2 f))(x) t^\ell + (\mathrm{Err}_t^{(m'-k,1)} J_k Q_t^2 f)(x).$$

As a result, we have

$$Q_t^1 Q_t^2 f(x) = \sum_{k=0}^{m'} \sum_{\ell=0}^{m'-k} (J_\ell Q_t^1 (J_k Q_t^2 f))(x) t^{k+\ell} + (\mathrm{Err}_t^{(m',1,2)} f)(x),$$

where

$$(5.1) \quad (\mathrm{Err}_t^{(m',1,2)} f)(x) = (Q_t^1 \mathrm{Err}_t^{(m',2)} f)(x) + \sum_{k=0}^{m'} (\mathrm{Err}_t^{(m'-k,1)} J_k Q_t^2 f)(x) t^k.$$

We obtain from the properties of the error terms that $\mathrm{Err}_t^{(m',1,2)} f \in C_{q'}$ for some $q' = q'(m,p) > q$.

(B): For $f \in C_p^{m''}$ with $m'' \geq 2(m'+1)$, we can derive for $k + \ell \leq m'$,

$$\|J_\ell Q_t^1 (J_k Q_t^2 f)\|_{C_{p+2(k+\ell)}^{m''-2(k+\ell)}} \leq K_1 \|J_k Q_t^2 f\|_{C_{p+2k}^{m''-2k}} \leq K_2 \|f\|_{C_p^{m''}}$$

and by (5.1),

$$\|\mathrm{Err}_t^{(m',1,2)} f\|_{C_{q'}} \leq K_3 \|\mathrm{Err}_t^{(m',2)} f\|_{C_q} + K_4 \|\mathrm{Err}_t^{(m',1)} J_0 Q_t^2 f\|_{C_{q'}}$$

$$+ K_5 \sum_{k=1}^{m'} \|J_k Q_t^2 f\|_{C_{p+2k}^{m''-2k}} t^{m'+1}$$

$$\leq \begin{cases} K t^{m'+1} \|f\|_{C_p^{m''}}, & \text{if } m' < m, \\ K t \delta_m(t) \|f\|_{C_p^{m''}}, & \text{if } m' = m. \end{cases}$$

Finally, the proof of (C) is straightforward. □

PROOF OF THEOREM 5.1. Using this lemma, we end the proof, calculating $J_{\leq m}$ for each numerical discretization scheme. For instance, in the case of N–V (b) [i.e., (2.10)], we obtain

$$J_{\leq 2}\bigg( \frac{1}{2} \prod_{i=0}^{d+1} Q_{i,t} + \frac{1}{2} \prod_{i=0}^{d+1} Q_{d+1-i,t} \bigg)$$

$$= \frac{1}{2} J_{\leq 2}\left(\prod_{i=0}^{d+1} J_{\leq 2}(Q_{it})\right) + \frac{1}{2} J_{\leq 2}\left(\prod_{i=0}^{d+1} J_{\leq 2}(Q_{d+1-i,t})\right)$$

$$= \frac{1}{2} J_{\leq 2}\left(\prod_{i=0}^{d+1}\left(\sum_{k=0}^{2} \frac{t^k}{k!} L_i^k\right)\right) + \frac{1}{2} J_{\leq 2}\left(\prod_{i=0}^{d+1}\left(\sum_{k=0}^{2} \frac{t^k}{k!} L_{d+1-i}^k\right)\right)$$

$$= \frac{1}{2}\left(I + t\sum_{i=1}^{d+1} L_i + \frac{t^2}{2}\sum_{i=1}^{d+1} L_i^2 + t^2 \sum_{i<j} L_i L_j\right)$$

$$+ \frac{1}{2}\left(I + t\sum_{i=1}^{d+1} L_i + \frac{t^2}{2}\sum_{i=1}^{d+1} L_i^2 + t^2 \sum_{i>j} L_i L_j\right)$$

$$= J_{\leq 2}(P_t). \qquad \square$$

Another idea to construct higher-order schemes is to use local Romberg extrapolation. In order to do this we need to weaken the assumption $\{\xi_i\} \subset [0,1]$. This is done in the next theorem.

THEOREM 5.2. *Let $m = 1$ or $2$. Assume (M) and $\mathcal{R}(2m, t^{2m})$ for $P_t$ and $Q_t^{[i]}$ ($i = 1, \ldots, \ell$) and (M$_P$) for $P_t$. Furthermore, we assume:*

(1) $J_{\leq 2m}(P_t - \sum_{i=1}^{\ell} \xi_i Q_t^{[i]}) = 0$ *for some real numbers $\{\xi_i\}_{i=1,\ldots,\ell}$ with $\sum_{i=1}^{\ell} \xi_i = 1$.*
(2) *There exists a constant $q = q(m, p) > 0$ such that for every $f \in C_p^{m'}$ with $m' \geq 2(m+1)$, $(P_t - Q_t^{[i]})f \in C_q^{m'-2(m+1)}$ and*

$$\sup_{t \in [0,T]} \|(P_t - Q_t^{[i]})f\|_{C_p^{m'-2(m+1)}} \leq C_T \|f\|_{C_q^{m'}} T^{m+1}.$$

*Then we have for any $f \in C_p^{4(m+1)}$,*

$$\left| P_T f(x) - \sum_{i=1}^{\ell} \xi_i (Q_{T/n}^{[i]})^n f(x) \right| \leq \frac{C(T, f, x)}{n^{2m}}.$$

PROOF.   We first remark that the operator $\sum_{i=1}^{\ell} \xi_i Q_t^{[i]}$ no longer satisfies the semigroup property, that is, $\sum_{i=1}^{\ell} \xi_i (Q_{T/n}^{[i]})^n \neq (\sum_{i=1}^{\ell} \xi_i Q_{T/n}^{[i]})^n$. Thus the proof is nontrivial.

Note that for $f \in C_p^{4(m+1)}$,

$$\mathcal{E} := P_T f(x) - \sum_{i=1}^{\ell} \xi_i (Q_{T/n}^{[i]})^n f(x) = \sum_{i=1}^{\ell} \xi_i (P_T - (Q_{T/n}^{[i]})^n) f(x).$$

Using the semigroup property of $P_t$ and $Q^{[i]}_{k/nT}$, we have

$$\mathcal{E} = \sum_{i=1}^{\ell} \xi_i \sum_{k=0}^{n-1} (Q^{[i]}_{T/n})^k (P_{T/n} - Q^{[i]}_{T/n}) P_{T-(k+1)/nT} f(x)$$

$$= \sum_{i=1}^{\ell} \xi_i \sum_{k=0}^{n-1} P_{kT/n} (P_{T/n} - Q^{[i]}_{T/n}) P_{T-(k+1)/nT} f(x)$$

$$+ \sum_{i=1}^{\ell} \xi_i \sum_{k=0}^{n-1} ((Q^{[i]}_{T/n})^k - P_{kT/n}) (P_{T/n} - Q^{[i]}_{T/n}) P_{T-(k+1)/nT} f(x).$$

We expand $(Q^{[i]}_{T/n})^k - P_{kT/n}$ again, to obtain

$$\mathcal{E} = \sum_{k=0}^{n-1} (P_{T/n})^k \left( P_{T/n} - \sum_{i=1}^{\ell} \xi_i Q^{[i]}_{T/n} \right) P_{T-(k+1)/nT} f(x)$$

$$+ \sum_{i=1}^{\ell} \xi_i \sum_{k=0}^{n-1} \sum_{l=0}^{k-1} (Q^{[i]}_{T/n})^l (Q^{[i]}_{T/n} - P_{T/n})$$

$$\times P_{T-(l+1)/nT} (P_{T/n} - Q^{[i]}_{T/n}) P_{T-(k+1)/nT} f(x).$$

By the assumption (1), we have

$$\left| \sum_{k=0}^{n-1} (P_{T/n})^k \left( P_{T/n} - \sum_{i=1}^{\ell} \xi_i Q^{[i]}_{T/n} \right) P_{T-(k+1)/nT} f(x) \right| \leq \frac{C_1(T,f,x)}{n^{2m}}.$$

Thus we end the proof by showing that

$$\left| \sum_{i=1}^{\ell} \xi_i \sum_{k=0}^{n-1} \sum_{l=0}^{k-1} (Q^{[i]}_{T/n})^l (Q^{[i]}_{T/n} - P_{T/n}) P_{T-(l+1)/nT} \right.$$

$$\left. \times (P_{T/n} - Q^{[i]}_{T/n}) P_{T-(k+1)/nT} f(x) \right|$$

$$\leq \frac{C_2(T,f,x)}{n^{2m}}.$$

Using here the assumption (2), we obtain

$$\|(Q^{[i]}_{T/n} - P_{T/n}) P_{T-(l+1)/nT} (P_{T/n} - Q^{[i]}_{T/n}) P_{T-(k+1)/nT} f\|_{C_{q'}}$$

$$\leq \frac{C(T)}{n^{m+1}} \|(P_{T/n} - Q^{[i]}_{T/n}) P_{T-(k+1)/nT} f\|_{C_q^{2(m+1)}}$$

$$\leq \frac{C'(T)}{n^{2(m+1)}} \|f\|_{C_p^{4(m+1)}}$$

and therefore

$$\left| \sum_{i=1}^{\ell} \xi_i \sum_{k=0}^{n-1} \sum_{l=0}^{k-1} (Q_{T/n}^{[i]})^l (Q_{T/n}^{[i]} - P_{T/n}) P_{T-(l+1)/nT} \right.$$

$$\left. \times (P_{T/n} - Q_{T/n}^{[i]}) P_{T-(k+1)/nT} f(x) \right|$$

$$\leq \sum_{k=0}^{n-1} \sum_{l=0}^{k-1} \frac{C_2(T, f, x)}{n^{2(m+1)}} \leq \frac{C_2(T, f, x)}{n^{2m}}.$$

This completes the proof.  □

EXAMPLE 5.2. It is known that the Ninomiya–Victoir scheme

$$\left( \frac{1}{2} e^{T/(2n)L_0} \prod_{i=1}^{d+1} e^{T/nL_i} e^{T/(2n)L_0} + \frac{1}{2} e^{T/(2n)L_0} \prod_{i=1}^{d+1} e^{T/nL_{d+2-i}} e^{T/(2n)L_0} \right)^n$$

is of order 2 [$m = 2$, $\delta_2(t) = t^2$ in Theorem 4.1]. By Theorem 5.2, the following modified Ninomiya–Victoir scheme

$$\frac{1}{2} \left( e^{T/(2n)L_0} \prod_{i=1}^{d+1} e^{T/nL_i} e^{T/(2n)L_0} \right)^n + \frac{1}{2} \left( e^{T/(2n)L_0} \prod_{i=1}^{d+1} e^{T/nL_{d+2-i}} e^{T/(2n)L_0} \right)^n$$

is also of order 2.

EXAMPLE 5.3. Fujiwara [6] gives a proof of a similar version of the above theorem and some examples of fourth and sixth order. We introduce the examples of fourth order:

$$\frac{4}{3} \left( \frac{1}{2} \left( \prod_{i=0}^{d+1} e^{t/2L_i} \right)^2 + \frac{1}{2} \left( \prod_{i=0}^{d+1} e^{t/2L_{d+1-i}} \right)^2 \right) - \frac{1}{3} \left( \frac{1}{2} \prod_{i=0}^{d+1} e^{tL_i} + \frac{1}{2} \prod_{i=0}^{d+1} e^{tL_{d+1-i}} \right).$$

In order to complete the approximation procedure through (quasi) Monte Carlo methods we need to find a stochastic characterization of the operators $Q_{i,t}$.

DEFINITION 5.4. Given a stochastic process $Y_t(x) \in \bigcap_{p \geq 1} L^p$, we say that $Y$ is the stochastic characterization of the linear operator $Q_t$ if $Q_t f(x) = E[f(Y_t(x))]$ for $f \in \bigcup_{p \geq 0} C_p$. In such as case we use the notation $Q_t \equiv Q_t^Y$.

REMARK 5.5. Given the operators $Q_t^{Z^i}$ ($i = 1, \ldots, \ell$) and the deterministic positive weights $\{\xi_i\}_{1 \leq i \leq \ell}$ with $\sum_{i=1}^{l} \xi_i = 1$. Let $U$ be a uniform random variable on $[0,1]$ independent of $(Z^i)_i$ and define $Z := \sum_{i=1}^{\ell} 1(\sum_{j=1}^{i-1} \xi_j \leq U < \sum_{j=1}^{i} \xi_j) Z^i$. Then

$$Q_t^Z f(x) \equiv E[f(Z_t(x))] = \sum_{i=1}^{\ell} \xi_i Q_t^{Z^i} f(x).$$

Therefore by Lemma 3.1 if $Q_t^{Z^i}$ satisfy $(\mathcal{M})$ and $\mathcal{R}(m, \delta_m)$ so does $Q_t^Z$. This property will be used repeatedly in what follows.

**6. Applications.** From this section on, we discuss the application of the previous approximation results to the case of solutions of the sde (1.1). From the results in the Appendix (see Corollary A.7), it is clear that the semigroup $P_t f(x) := E[f(X_t(x))]$ satisfies the hypotheses $(\mathcal{M})$ and $\mathcal{R}(m, \delta_m)$. We define various approximations generated via a stochastic process $\bar{X}_i$ with corresponding operator $Q_t^{\bar{X}_i}$ $(i = 0, 1, \ldots, d+1)$.

Due to the previous results and in particular, Theorem 5.1, we see that is enough to verify local conditions on the approximation operators to conclude global properties of approximation. In particular, we only need to verify that the operator associated with $\bar{X}_i$ (the approximation to the coordinate process) satisfies $(\mathcal{M})$ and $\mathcal{R}(m, \delta_m)$ and $J_{\leq m}(Q_t^{\bar{X}_i}) = I + \sum_{j=1}^{m} \frac{t^j}{j!} L_i^j$ for some $m \geq 2$ for $L_i$ given by (2.5). This is the goal in most of the applications in this section.

Recall that the stochastic differential equation to be approximated is

$$X_t(x) = x + \sum_{i=0}^{d} \int_0^t V_i(X_{s-}(x)) \circ dB_s^i + \int_0^t h(X_{s-}(x)) \, dY_s.$$

In each of the following sections we consider different approximation processes for the coordinate processes $X_{i,t}$. In each section, the notation for the approximating process is always $\bar{X}_{i,t}$. We hope that this does not raise confusion as the framework in each section is clear.

6.1. *Continuous diffusion component.*

(a) *Explicit solution.* Let $V : \mathbf{R}^N \to \mathbf{R}^N$ be a smooth function satisfying the linear growth condition $|V(x)| \leq C(1 + |x|)$. The exponential map is defined as $\exp(V)x = z_1(x)$ where $z$ denotes the solution of the ordinary differential equation

(6.1) $$\frac{dz_t(x)}{dt} = V(z_t(x)), \qquad z_0(x) = x.$$

The solution of the coordinate sde is obtained in the following proposition. The proof follows from Itô's formula.

PROPOSITION 6.1. *For $i = 0, 1, \ldots, d$, the stochastic differential equation*

$$X_{i,t}(x) = x + \int_0^t V_i(X_{i,s}(x)) \circ dB_s^i \tag{6.2}$$

*has a unique solution given by*

$$X_{i,t}(x) = \exp(B_t^i V_i)x.$$

$X_{i,t}(x)$ is called the $i$th coordinate process and its semigroup is denoted by $Q_t^i$. This is a trivial example of the approximation of $e^{tL_i}$, $i = 0, 1, \ldots, d$, satisfying $(\mathcal{M})$ and $\mathcal{R}(m, t^m)$. However, sometimes it is not easy to obtain the closed-form solution to the ODE (6.1). In those cases, we shall approximate $\exp(tV)x$. Here we will do this with the Taylor expansion first and then the Runge–Kutta methods denoted by $b_m$ and $c_m$, respectively.

(b) *Taylor expansion.* We first prove the following lemmas which help us to find the rate of convergence of the scheme to be defined later. The following lemma follows easily from Gronwall's lemma.

LEMMA 6.2. *Let $V$ be a smooth function which satisfies the linear growth condition. Then $|\exp(tV)x| \leq Ce^{K|t|}(1 + |x|)$ for $t \in \mathbf{R}$, $x \in \mathbf{R}^N$.*

From now on we denote by $e_j : \mathbf{R}^N \to \mathbf{R}$, the coordinate function $e_j(x) = x_j$ for $j = 1, \ldots, N$. Furthermore, we also denote by $V$ the vector field operator defined from $V$.

LEMMA 6.3. *Let $f \in C_p^{m+1}$. Then we have for $i = 0, 1, \ldots, d$,*

$$f(\exp(tV_i)x) = \sum_{k=0}^m \frac{t^k}{k!} V_i^k f(x) + \int_0^t \frac{(t-u)^m}{m!} V_i^{m+1} f(\exp(uV_i)x) \, du \tag{6.3}$$

*and*

$$\left| \int_0^t \frac{(t-u)^m}{m!} V_i^{m+1} f(\exp(uV_i)x) \, du \right| \leq C_m \|f\|_{C_p^{m+1}} e^{K|t|} (1 + |x|^{p+m+1}) |t|^{m+1}$$

*for all $t \in \mathbf{R}$.*

PROOF. Assertion (6.3) follows application of Taylor expansion to the function $f(\exp(tV)x)$ around $t = 0$. Next, as $|V_i^{m+1} f(x)| \leq C(1 + |x|^{p+m+1})$, we obtain from Lemma 6.2,

$$\left| \int_0^t \frac{(t-u)^m}{m!} V_i^{m+1} f(\exp(uV)x) \, du \right|$$

$$\leq C_m \|f\|_{C_p^{m+1}} \int_0^{|t|} |t|^m C e^{K|u|} (1 + |x|^{p+m+1}) \, du$$

$$\leq C_m' \|f\|_{C_p^{m+1}} e^{K|t|} (1 + |x|^{p+m+1}) |t|^{m+1}.$$

□

Based on this lemma, we define the approximation to the solution of the coordinate equation (6.2) as follows

$$b_m^j(t,V)x = \sum_{k=0}^{m} \frac{t^k}{k!}(V^k e_j)(x), \qquad j=1,\ldots,N.$$

Define

$$\bar{X}_{i,t}(x) = b_{2m+1}(B_t^i, V_i)x \qquad \text{for } i=0,\ldots,d.$$

Then we have the following approximation result.

PROPOSITION 6.4.  (i) *For every* $p \geq 1$,

$$\|X_{i,t}(x) - \bar{X}_{i,t}(x)\|_{L^p} \leq C(p,m,T)(1+|x|^{2(m+1)})t^{m+1}.$$

(ii) *Let* $f \in C_p^1$. *Then we have*

$$E[|f(X_{i,t}(x)) - f(\bar{X}_{i,t}(x))|] \leq C(m,T)\|f\|_{C_p^1}(1+|x|^{p+2(m+1)})t^{m+1}.$$

PROOF.  (i) Apply Proposition 6.1 and Lemma 6.3 with $f = e_i$. Then we have

$$\|X_{i,t}(x) - \bar{X}_{i,t}(x)\|_{L^p} \leq E[|C_m e^{K|B_t|}(1+|x|^{2(m+1)})|B_t|^{2(m+1)}|^p]^{1/p}$$

$$\leq C(1+|x|^{2(m+1)})t^{m+1}$$

for some constant $C = C(p,m,T)$.

(ii) We first apply the mean value theorem to obtain

$$E[|f(X_{i,t}(x)) - f(\bar{X}_{i,t}(x))|]$$

$$\leq \|f\|_{C_p^1}\|1 + |\theta X_{i,t}(x) + (1-\theta)\bar{X}_{i,t}(x)|^p\|_{L^2}\|X_{i,t}(x) - \bar{X}_{i,t}(x)\|_{L^2}$$

$$\leq C\|f\|_{C_p^1}\|1 + |X_{i,t}(x)|^p + |\bar{X}_{i,t}(x)|^p\|_{L^2}(1+|x|^{2(m+1)})t^{m+1}.$$

We see by Lemma 6.2 that

$$\sup_{t\in[0,T]} \|1 + |X_{i,t}(x)|^p + |\bar{X}_{i,t}(x)|^p\|_{L^2} \leq C'(1+|x|^p)$$

from which the proof follows.  □

As a result of this proposition we can see that $\mathcal{R}(m,t^m)$ holds for the operators associated with $b_m(t,V_0)x$ and $b_{2m+1}(B_t^i,V_i)x$, $1 \leq i \leq d$. Indeed, we have for $m' \leq m$,

$$E[f(\bar{X}_{i,t}(x))] = E[f(X_{i,t}(x))] + E[f(\bar{X}_{i,t}(x)) - f(X_{i,t}(x))]$$

$$= \sum_{k=0}^{m'} \frac{t^k}{k!} L_i^k f(x) + (E_t^{m'} f)(x),$$

where

$$(E_t^{m'} f)(x) := (\mathrm{Err}_t^{(m')} f)(x) + E[f(\bar{X}_{i,t}(x)) - f(X_{i,t}(x))]$$

and $(\mathrm{Err}_t^{(m')} f)(x)$ is defined through the residue appearing in Proposition A.6, using $L_i$ and $Q_i$ instead of $L$ and $P$. Furthermore, using Proposition 6.4(ii), we have that the error term $E_t^{m'}$ satisfies (B) in Assumption $\mathcal{R}(m, t^m)$.

It remains to prove that $(\mathcal{M})$ holds for $\bar{X}_{i,t}(x)$. For the proof, we need an additional growth condition for the vector field $V_i$.

PROPOSITION 6.5. *Assume that $(V_i^k e_j)$ $(2 \le k \le m,\ 0 \le i \le d,\ 1 \le j \le N)$ satisfies the linear growth condition then $(\mathcal{M})$ holds for $\bar{X}_{i,t}(x)$, $i = 0, \ldots, d$.*

PROOF. The Assumption $(\mathcal{M}_0)$ follows from the smoothness and the linear growth property of $V_i^k e_j$. We only prove the moment condition (3.1) for $\bar{X}_{i,t}(x)$ $i = 1, \ldots, d$. Consider the multiplication ($p \in \mathbf{N}$)

$$\left| \sum_{k=0}^{m} \frac{(B_t^i)^k}{k!} (V_i^k e_j)(x) \right|^{2p} = \left| x + B_t^i V_i(x) + \sum_{k=2}^{m} \frac{(B_t^i)^k}{k!} (V_i^k e_j)(x) \right|^{2p}.$$

Taking into account that $E[(B_t^i)^{2k+1}] = 0$, $k \in \mathbf{N}$. Then by the assumption, we obtain the result. □

Therefore we obtain the main result.

THEOREM 6.1. *Assume that $(V_i^k e_j)$ $(2 \le k \le m,\ 0 \le i \le d,\ 1 \le j \le N)$ satisfies the linear growth condition. Let $\bar{X}_{i,t}(x)$ be defined by*

$$\bar{X}_{i,t}(x) = b_{2m+1}(B_t^i, V_i)x = \sum_{k=0}^{2m+1} \frac{1}{k!} (V_i^k I)(x) \int_{0 < t_1 < \cdots < t_k < t} 1 \circ dB_{t_1}^i \circ \cdots \circ dB_{t_k}^i.$$

*Denote by $Q_t^{\bar{X}_i}$ the semigroup associated with $\bar{X}_{i,t}(x)$. Then $Q_t^{\bar{X}_i}$ satisfies $(\mathcal{M})$ and $\mathcal{R}(m, t^m)$. Furthermore $J_{\le m}(Q_t^{\bar{X}_i}) = I + \sum_{j=1}^{m} \frac{t^j}{j!} L_i^j$.*

(c) *Runge–Kutta methods.* We say here that $c_m$ is an $s$-stage explicit Runge-Kutta method of order $m$ for the ODE (6.1) if it can be written in the form

(6.4)  $$c_m(t, V)x = x + t \sum_{i=1}^{s} \beta_i k_i(t, V)x,$$

where $k_i(t, V)x$ defined inductively by

$$k_1(t, V)x = V(x),$$

$$k_i(t, V)x = V\left(x + t\sum_{j=1}^{i-1} \alpha_{i,j} k_j(t, V)x\right), \qquad 2 \leq i \leq s,$$

and satisfies

$$|\exp(tV)x - c_m(t, V)x| \leq C_m e^{K|t|}(1 + |x|^{m+1})|t|^{m+1}$$

for some constants $((\beta_i, \alpha_{i,j})_{1 \leq j < i \leq s})$. Runge–Kutta formulas of order less than or equal to 7 are well known. For details, see, for example, Butcher [4].

The following proposition can be shown by the same argument as in the proof of Proposition 6.4.

PROPOSITION 6.6 (Stochastic Runge–Kutta). (i) *For every $p \geq 1$,*

(6.5) $\quad \|X_{i,t}(x) - c_{2m+1}(B_t^i, V_i)x\|_{L^p} \leq C(p, m, T)(1 + |x|^{2(m+1)})t^{m+1}.$

(ii) *Let $f \in C_p^1$. Then we have*

(6.6)
$$E[|f(X_{i,t}(x)) - f(c_{2m+1}(B_t^i, V_i)x)|]$$
$$\leq C(m, T)\|f\|_{C_p^1}(1 + |x|^{2(m+1)})t^{m+1}.$$

Next we show that $(\mathcal{M})$ still holds for the Runge–Kutta schemes.

PROPOSITION 6.7. *$(\mathcal{M})$ holds for $c_m(B_t^i, V_i)x$, $i = 0, \ldots, d$.*

PROOF. We first note that for every $1 \leq j \leq s$, there exists a function of the form $p_j = \sum_{k=0}^{j-1} a_{jk}|t|^k$ such that

$$|k_j(t, V)x| \leq p_j(t)(1 + |x|).$$

Assumption $(\mathcal{M}_0)$ follows from the smoothness and the linear growth property of $V_i$. We now prove (3.1). In the case $i = 0$, this is obvious by definition and the inequality (6.1). In the case $1 \leq i \leq d$, observe that

$$c_m(t, V)x = x + t\sum_{l=1}^{s} \beta_l V(x) + t\sum_{l=2}^{s} \beta_l \int_0^1 \frac{d}{d\theta} V\left(x + \theta t \sum_{j=1}^{l-1} \alpha_{l,j} k_j(t, V)x\right) d\theta$$

$$=: x + t\sum_{l=1}^{s} \beta_l V(x) + D_m(t, V)x.$$

Expanding multiplications and taking expectations, as in Proposition 6.5, we can show that the terms containing odd powers of $B_t^i$ have expectation 0. Finally, we obtain from the boundedness of $\partial V_i$ that

$$|D_m(B_t^i, V_i)x| \leq p(B_t^i)(1 + |x|),$$

where $p = p(t)$ is of the form $\sum_{k=2}^{s} a_k |t|^k$. Using this, we conclude the proof. □

Consequently, as in the Taylor scheme, $\mathcal{R}(m, t^m)$ and $(\mathcal{M})$ hold for the operators associated with $c_m(t, V_0)x$ and $c_{2m+1}(B_t^i, V_i)x$, $1 \le i \le d$. For more on this method, we refer the reader to [12].

(d) *Minor extension.* In the previous approximation, the assumption that $B_t \sim N(0, I_d)$ can be weakened. In fact, we can use $\sqrt{t}Z$ instead of $B_t$ where $(Z^i)_{i=1}^d$ are independent and

$$P(Z^i = \pm\sqrt{3}) = \tfrac{1}{6}, \qquad P(Z^i = 0) = \tfrac{2}{3}$$

for each $i = 1, \ldots, d$.

PROPOSITION 6.8. *Let $B_t$ be a one-dimensional Brownian motion and $Z$ be a $\mathbf{R}$-valued random variable such that for all $0 \le k \le m$,*

$$E[(Z)^k] = E[(B_1)^k]$$

*and*

$$E[\exp(c|Z|)] < \infty$$

*for any $c > 0$. Then, for every $f \in C_p^{m+1}$,*

$$|E[f(\exp(B_t V)x)] - E[f(c_m(\sqrt{t}Z, V)x)]| \le C(m, T)(1 + |x|^{p+m+1})t^{(m+1)/2}.$$

6.2. *Compound Poisson case.* Suppose that $Y_t$ is a compound Poisson process. That is,

$$Y_t = \sum_{i=1}^{N_t} J_i,$$

where $(N_t)$ is a Poisson process with intensity $\lambda$ and $(J_i)$ are i.i.d. $\mathbf{R}^d$-valued random variables independent of $(N_t)$ with $J_i \in \bigcap_{p \ge 1} L^p$.

In this case $Y_t$ is a Lévy process with generator of the form

$$\int_{\mathbf{R}_0^d} (f(x+y) - f(x))\nu(dy),$$

where $\tau \equiv 0$, $b = 0$, $\nu(\mathbf{R}_0^d) = \lambda < \infty$ and $\nu(dy) = \lambda P(J_1 \in dy)$.

Then in this case

(6.7) $$X_t^{d+1}(x) = x + \int_0^t h(X_{s-}^{d+1}(x)) \, dY_s, \qquad t \in [0, T],$$

which can be solved explicitly. Indeed, let $(G_i(x))$ be defined by recursively

$$G_0 = x,$$
$$G_i = G_{i-1} + h(G_{i-1})J_i.$$

Then the solution can be written as $X_t^{d+1}(x) = G_{N_t}(x)$. Define for fixed $M \in \mathbf{N}$, the approximation process $\bar{X}_{d+1,t} = G_{N_t \wedge M}(x)$. This approximation requires the simulation of at most $M$ jumps. In fact, the rate of convergence is fast as the following result shows (see also Mordecki et al. [14]).

THEOREM 6.2. *Let $M \in \mathbf{N}$. Then the process $G_{N_t \wedge M}(x)$ satisfies $(\mathcal{M})$ and $\mathcal{R}(M, t^{M-\kappa})$ for arbitrary small $\kappa > 0$. Furthermore $J_{\leq M}(Q_t^{\bar{X}_{d+1}}) = I + \sum_{j=1}^m \frac{t^j}{j!} L_{d+1}^j$.*

PROOF. Note that for $f \in C_p$

$$Q_t^{\bar{X}_{d+1}} f(x) - Q_t^{d+1} f(x) = E[f(G_{N_t \wedge M}(x))] - E[f(G_{N_t}(x))]$$
$$= E[(f(G_{N_t \wedge M}(x)) - f(G_{N_t}(x)))1_{\{T_{M+1} \leq t\}}],$$

where $T_M := \inf\{t > 0 : N_t = M\}$. By the Hölder inequality,

$$|Q_t^{\bar{X}_{d+1}} f(x) - Q_t^{d+1} f(x)|$$
$$\leq 2E\left[\sup_{0 \leq t \leq T} |f(G_{N_t}(x))|^{\gamma/(\gamma-1)}\right]^{(\gamma-1)/\gamma} P(T_{M+1} \leq t)^{1/\gamma}$$
$$= 2E\left[\sup_{0 \leq t \leq T} |f(G_{N_t}(x))|^{\gamma/(\gamma-1)}\right]^{(\gamma-1)/\gamma} \left(\int_0^t \frac{(\lambda s)^M}{M!} \lambda e^{-\lambda s} ds\right)^{1/\gamma}$$
$$\leq C(\gamma, T) \|f\|_{C_p} (1 + |x|^p)(t\lambda^{-1})^{(M+1)/\gamma}.$$

Take sufficiently small $\gamma > 1$, then $\mathcal{R}(M, t^{M-\kappa})$ holds for $Q_t^{\bar{X}_{d+1}}$ where $\kappa := (1 - 1/\gamma)(M + 1) > 0$. Finally, we show $(\mathcal{M})$. Let $f_p(x) = |x|^{2p}$ ($p \in \mathbf{N}$) and $\gamma < M$. Then using the above calculation and Corollary A.7, we have

$$Q_t^{\bar{X}_{d+1}} f_p(x) = Q_t^{d+1} f_p(x) + (Q_t^{\bar{X}_{d+1}} f_p(x) - Q_t^{d+1} f_p(x))$$
$$\leq (1 + K_1 t) f_p(x) + K_2 t + |Q_t^{\bar{X}_{d+1}} f_p(x) - Q_t^{d+1} f_p(x)|$$
$$\leq (1 + K_3 t) f_p(x) + K_4 t. \qquad \square$$

6.3. *Infinite activity case.* In this section, we consider the SDE (2.3) under the conditions $\nu(\mathbf{R}_0^d) = \infty$. Without loss of generality, we assume that $c \equiv 0$.

(a) *Ignoring small jumps.* Define for $\varepsilon > 0$ the finite activity (i.e., drift + compound Poisson) Lévy process $(Y_t^\varepsilon)$ with Lévy triplet $(b, 0, \nu^\varepsilon)$ where the Lévy measure is defined by

$$\nu^\varepsilon(E) := \nu(E \cap \{y : |y| > \varepsilon\}), \qquad E \in \mathcal{B}(\mathbf{R}_0^d). \tag{6.8}$$

Consider the approximate coordinate SDE

$$\bar{X}_{d+1,t}(x) = x + \int_0^t h(\bar{X}_{d+1,s-}(x)) \, dY_s^\varepsilon,$$

result shows (see also Mordecki et al. [14]) whose generator is

$$L_{d+1}^{1,\varepsilon} f(x) = \nabla f(x) h(x) b + \int (f(x + h(x)y) - f(x) - \nabla f(x) h(x) \tau(y)) \nu^\varepsilon(dy).$$

Now we derive the error estimate for $\bar{X}_{d+1,t}$.

THEOREM 6.3. *Assume that $0 < \varepsilon \equiv \varepsilon(t) \leq 1$ is chosen as to satisfy that $\sigma^2(\varepsilon) := \int_{|y| \leq \varepsilon} |y|^2 \nu(dy) \leq t^{M+1}$. Then we have that $Q_t^{\bar{X}_{d+1}}$ satisfies $(\mathcal{M})$ and $\mathcal{R}(M, t^M)$. Furthermore $J_{\leq M}(Q_t^{\bar{X}_{d+1}}) = I + \sum_{j=1}^m \frac{t^j}{j!} L_{d+1}^j$.*

PROOF. First, we remark that condition $(\mathcal{M}_0)$ follows from Proposition 5.2 in [7]. We start by noting that from Proposition A.6, we have

$$Q_t^{d+1} f(x) - Q_t^{\bar{X}_{d+1}} f(x)$$

$$= \sum_{k=1}^M \frac{t^k}{k!} ((L_{d+1})^k - (L_{d+1}^{1,\varepsilon})^k) f(x) \tag{6.9}$$

$$+ \int_0^t \frac{(t-u)^M}{M!} (Q_u^{d+1} (L_{d+1})^{M+1} - Q_u^{\bar{X}_{d+1}} (L_{d+1}^{1,\varepsilon})^{M+1}) f(x) \, du.$$

Therefore the proof is achieved if we prove that

$$|(L_{d+1} - L_{d+1}^{1,\varepsilon}) f(x)| \leq C \|f\|_{C_p^2} (1 + |x|^{p+2}) t^{M+1}.$$

For the proof, we change here the representation of the Lévy triplets of $Y_t$ and $Y_t^\varepsilon$ as follows:

$$(b, 0, \nu), \tau \Rightarrow (b_\varepsilon, 0, \nu), \tau_\varepsilon,$$
$$(b, 0, \nu^\varepsilon), \tau \Rightarrow (b_\varepsilon, 0, \nu^\varepsilon), \tau_\varepsilon,$$

where $\tau_\varepsilon(y) = y 1_{\{|y| \leq \varepsilon\}}$. Then

$$|(L_{d+1} - L_{d+1}^{1,\varepsilon}) f(x)|$$

$$\leq \left| \int \nabla f(x) h(x) (y - \tau_\varepsilon(y)) (\nu(dy) - \nu^\varepsilon(dy)) \right| \tag{6.10}$$

$$+ \left| \int \int_0^1 (1 - \theta) \frac{d^2}{d\theta^2} f(x + \theta h(x) y) \, d\theta (\nu(dy) - \nu^\varepsilon(dy)) \right|.$$

We first obtain that for $\varepsilon > 0$,
$$\int (y - \tau_\varepsilon(y))(\nu(dy) - \nu^\varepsilon(dy)) = 0$$
since the support of the measure $(\nu - \nu^\varepsilon)$ is $\{|y| \leq \varepsilon\}$. Now we consider the second term of (6.10). We can immediately show that due to the polynomial growth property for $f$,
$$\left| \int\int_0^1 \frac{d^2}{d\theta^2} f(x + \theta h(x)y) \, d\theta (\nu(dy) - \nu^\varepsilon(dy)) \right| \leq C\|f\|_{C_p^2}(1 + |x|^{p+2})\sigma^2(\varepsilon)$$
and hence as $\sigma^2(\varepsilon) \leq t^{M+1}$, one obtains that $J_{\leq M}(Q_t^{\bar{X}_{d+1}}) = I + \sum_{j=1}^m \frac{t^j}{j!} L_{d+1}^j$ and that $Q_t^{\bar{X}_{d+1}}$ satisfies ($\mathcal{M}$) and $\mathcal{R}(M, t^M)$ follows as in the proof of Proposition 6.2. □

Using Theorem 5.1, we can incorporate the above approximating process $\bar{X}_{d+1,t}$ to the whole approximation method. This will require to first simulate the jump times of the approximating Lévy process $Y^\varepsilon$ and then solving ODEs between these times. If the task is time consuming one can also separate the jump component from the drift component as indicated by Theorem 5.1 (see also Section 6.4). The right size of $\varepsilon$ is determined by the condition $\sigma^2(\varepsilon) \leq t^{M+1}$.

(b) *Approximation of small jumps.* We apply here the Asmussen–Rosiński's approximation for small jumps of Lévy processes. The idea is that the small jumps ignored in (6.8) are close to a Brownian motion with small variance $\sigma^2(\varepsilon)$ (see details in [2]).

Consider the new approximate SDE

(6.11) $\quad \bar{X}_{d+1,t}(x) = x + \int_0^t h(\bar{X}_{d+1,s}(x)) \Sigma_\varepsilon^{1/2} \, dW_s + \int_0^t h(\bar{X}_{d+1,s-}(x)) \, dY_s^\varepsilon,$

where $W_t$ is a new $d$-dimensional Brownian motion independent of $B_t$ and $Y_t^\varepsilon$, and $\Sigma_\varepsilon$ is the symmetric and semi-positive definite $d \times d$ matrix defined as

(6.12) $$\Sigma_\varepsilon = \int_{|y| \leq \varepsilon} yy^* \nu(dy).$$

We remark that $\Sigma_\varepsilon$ is of the form $A\Lambda A^*$, where $A$ is an orthogonal matrix and $\Lambda$ is the diagonal matrix with entries $\lambda_1, \ldots, \lambda_d \geq 0$ (eigenvalues). Thus we use the notation $\Sigma_\varepsilon^{1/2} = A\Lambda^{1/2}$. Since the above SDE is also driven by a jump-diffusion process, we can also simulate it using the second-order discretization schemes in Theorem 5.1.

Now we prove that rate of convergence in this case is faster than in the case that we ignore completely the small jumps (see Theorem 6.3).

THEOREM 6.4. *Assume that $0 < \varepsilon \equiv \varepsilon(t) \leq 1$ is chosen as to satisfy that $\int_{|y|\leq \varepsilon} |y|^3 \nu(dy) \leq t^{M+1}$. Then we have that $Q_t^{\bar{X}_{d+1}}$ satisfies $(\mathcal{M})$ and $\mathcal{R}(M, t^M)$. Furthermore $J_{\leq M}(Q_t^{\bar{X}_{d+1}}) = I + \sum_{j=1}^m \frac{t^j}{j!} L_{d+1}^j$.*

PROOF. As before, condition $(\mathcal{M}_0)$ follows from Proposition 5.2 in [7]. The SDE $\bar{X}_{d+1,t}$ corresponds to the generator

$$L_{d+1}^{2,\varepsilon} f(x) := \nabla f(x) h(x) b + \frac{1}{2} \sum_{k,l} \partial_{k,l} f(x) (h(x) \Sigma_\varepsilon h^*(x))_{k,l}$$

$$+ \int (f(x + h(x)y) - f(x) - \nabla f(x) h(x) \tau(y)) \nu^\varepsilon(dy).$$

Using this representation, we have for $f \in C_p^3$,

$$(L_{d+1} - L_{d+1}^{2,\varepsilon}) f(x) = \iint_0^1 (1-\theta) \frac{d^2}{d\theta^2} f(x + \theta h(x) y) \, d\theta (\nu(dy) - \nu^\varepsilon(dy))$$

$$- \frac{1}{2} \sum_{k,l} \partial_{k,l} f(x) (h(x) \Sigma_\varepsilon h^*(x))_{k,l}$$

$$= \iint_0^1 \frac{(1-\theta)^2}{2} \frac{d^3}{d\theta^3} f(x + \theta h(x) y) \, d\theta (\nu(dy) - \nu^\varepsilon(dy)).$$

Hence we finish the proof as in the proof of Theorems 6.2 and 6.3. □

If we put all the pieces together, we have the following final result. Here $B_t^{ij}$ denote $i = 1, \ldots, d$, $j = 1, \ldots, 2n$ denote $2nd$ independent standard Brownian motions and $B_t^{0j} \equiv t$.

THEOREM 6.5. *Assume that $V_0$, $V$ and $h$ are infinitely differentiable functions with bounded derivatives with $\int_{\mathbf{R}_0^d} (1 \wedge |y|^p) \nu(dy) < \infty$ for all $p \in \mathbf{N}$. Define $\varepsilon \equiv \varepsilon(T, n)$ so that $\int_{|y|\leq \varepsilon} |y|^3 \nu(dy) \leq (\frac{T}{n})^3$. Let $\bar{X}_{i,t}^j(x) = c_5(B_t^{ij}, V_i)x$, $i = 0, \ldots, d$, $j = 1, \ldots, 2n$, $2n$ copies of the Runge–Kutta method of order 2 as defined in (6.4) and $\bar{X}_{d+1,t}^j(x)$, $j = 1, \ldots, 2n$, independent copies of the approximation defined in (6.11). Then the following schemes, $X_T^{(n)} = Y_n^n \circ Y_n^{n-1} \circ \cdots \circ Y_n^1(x)$, are second-order discretization schemes:*

*N–V (a): $Y_n^j(x) = U_j \bar{X}_{0,T/(2n)}^j \circ \bar{X}_{1,T/n}^j \circ \cdots \circ \bar{X}_{d+1,T/n}^j \circ \bar{X}_{0,T/(2n)}^j(x) + (1 - U_j) \bar{X}_{0,T/(2n)}^j \circ \bar{X}_{d+1,T/n}^j \circ \cdots \circ \bar{X}_{1,T/n}^j \circ \bar{X}_{0,T/(2n)}^j(x)$ where $U_j$ is a Bernoulli r.v. with $P(U_j = 1) = 1/2$, independent of everything else.*

*N–V (b): $Y_n^j(x) = U_j \bar{X}_{d+1,T/n}^j \circ \cdots \circ \bar{X}_{0,T/n}^j(x) + (1 - U_j) \bar{X}_{0,T/n}^j \circ \cdots \circ \bar{X}_{d+1,T/n}^j(x)$ where $U_j$ is a Bernoulli r.v. with $P(U_j = 1) = 1/2$, independent of everything else.*

Splitting: $Y_n^j(x) = \bar{X}_{0,T/(2n)}^j \circ \cdots \circ \bar{X}_{d,T/(2n)}^j \circ \bar{X}_{d+1,T/n}^j \circ \bar{X}_{d,T/(2n)}^{n+j} \circ \cdots \circ \bar{X}_{0,T/(2n)}^{n+j}(x)$.

One can also write a similar result for higher-order schemes using Theorem 5.2.

6.4. *Limiting the number of jumps per interval for approximations of infinite activity Lévy driven SDEs.* In the previous two approximations although $\varepsilon \in (0,1)$ may be relatively large compared with the interval size $T/n$, one still faces the possibility of having many jumps in the interval $[0,T]$. Therefore we introduce the idea used in Section 6.2. That is, we propose another approximation that restricts the numbers of possible jumps to at most $n$. Throughout this section we assume that $\int_{|y|<1} |y|\nu(dy) < \infty$ and without loss of generality, we assume that $\tau(y) = y1_{|y|<1}$.

Then we decompose the operator

$$L_{d+1} = L_{d+1}^1 + L_{d+1}^2 + L_{d+1}^3,$$

$$L_{d+1}^1 f(x) := \nabla f(x) h(x) \left( b - \int_{\varepsilon < |y| \leq 1} \tau(y) \nu(dy) \right),$$

$$L_{d+1}^2 f(x) := \int_{|y| \leq \varepsilon} (f(x + h(x)y) - f(x) - \nabla f(x) h(x) \tau(y)) \nu(dy),$$

$$L_{d+1}^3 f(x) := \int_{\varepsilon < |y|} f(x + h(x)y) - f(x) \nu(dy).$$

The operator $L_{d+1}^1$ can be easily approximated using any Runge–Kutta method for the ordinary differential equation

$$X_{d+1,t}^1 = x + \left( b - \int_{\varepsilon < |y| \leq 1} \tau(y) \nu(dy) \right) \int_0^t h(X_{d+1,s}^1) \, ds.$$

We denote by $\bar{X}_{d+1,t}^1$, the Euler scheme associated with this ordinary differential equation. Therefore we only need to approximate $L_{d+1}^2$ and $L_{d+1}^3$.

Let $l : \mathbf{R}^d \to \mathbf{R}_+$ be a localization function that may be used for importance sampling of the Lévy measure. Let $F_\varepsilon^l(dy) = \lambda_\varepsilon^{-1} l(y) 1_{|y| \leq \varepsilon} \nu(dy)$ with $\lambda_\varepsilon = \int_{|y| \leq \varepsilon} l(y) \nu(dy)$. Let $Y_\varepsilon \sim F_\varepsilon$. Define $\bar{X}_t^{2,\varepsilon}(x) = x + h(x) W_t \sqrt{\lambda_\varepsilon}$, where $W$ is a $d$-dimensional Wiener process with covariance matrix given by $\Sigma_{ij} = l(Y^\varepsilon)^{-1} Y_i^\varepsilon Y_j^\varepsilon$ which is independent of everything else.

First we prove that $\bar{X}_t^{2,\varepsilon}(x)$ satisfies Assumption $(\mathcal{M})$.

LEMMA 6.9. *Assume that for $p \geq 2$, $\sup_{\varepsilon \in (0,1]} \int_{|y| \leq \varepsilon} |y|^p l(y)^{-(p-2)/2} \times \nu(dy) < \infty$, then Assumption $(\mathcal{M})$ is satisfied with*

$$E[|\bar{X}_{d+1}^{2,\varepsilon}(x)|^p] \leq (1 + Kt)|x|^p + K't.$$

PROOF. Let $f(x) = |x|^p$, $p \geq 2$. Using Itô's formula for $p \neq 3$ and an approximative argument in the case $p = 3$ (as in the proof of the Meyer–Itô formula) one obtains that

$$\text{(6.13)} \quad E[f(\bar{X}_t^{2,\varepsilon}(x))] - f(x)$$

$$\text{(6.14)} \quad = \frac{p}{2}\lambda_\varepsilon E\left[l(Y^\varepsilon)^{-1}\int_0^t \left(\frac{p}{2}-1\right)|\bar{X}_s^{2,\varepsilon}(x)|^{p-4}\langle h(x)Y^\varepsilon, \bar{X}_s^{2,\varepsilon}(x)\rangle^2 \right.$$

$$\left. + |\bar{X}_s^{2,\varepsilon}(x)|^{p-2}|h(x)Y^\varepsilon|^2\, ds\right].$$

We use the Lipschitz property of $h$ to obtain that

$$|\bar{X}_s^{2,\varepsilon}(x)| = |x + h(x)W_s\sqrt{\lambda_\varepsilon}|$$
$$\leq (1 + C|W_s|\sqrt{\lambda_\varepsilon})(1 + |x|).$$

Then, we have

$$|E[f(\bar{X}_t^{2,\varepsilon}(x))] - f(x)|$$
$$\leq C_p t(1 + |x|^p)\int_{|y|<\varepsilon} |y|^2(1 + (|y|^2 l(y))^{-1}\lambda_\varepsilon t)^{(p-2)/2})\nu(dy). \qquad \square$$

LEMMA 6.10.  *Assume that for $p \geq 2$,*

$$M_p = \sup_{\varepsilon \in (0,1]}\int_{|y|\leq\varepsilon} |y|^4 l(y)^{-1}(1 + (|y|^2 l(y))^{-1}\lambda_\varepsilon t)^{(p-2)/2})\nu(dy) < \infty$$

*and $\int_{|y|\leq\varepsilon} |y|^3\nu(dy) \leq Ct$ then*

$$|E[f(\bar{X}_t^{2,\varepsilon})] - f(x) - tL_{d+1}^2 f(x)| \leq C(p)\|f\|_{C_p^4}(1 + |x|^{p+4})t^2.$$

*That is, $\bar{X}_t^{2,\varepsilon}(x)$ satisfies Assumption $\mathcal{R}(2, t^2)$.*

PROOF. Let $f \in C_p^4$ then applying Itô's formula, one gets

$$E[f(\bar{X}_t^{2,\varepsilon})] = f(x) + \frac{\lambda_\varepsilon}{2}E\left[\int_0^t \sum_{i,j,k,l} \partial_{ij}f(\bar{X}_s^{2,\varepsilon})h_{ik}h_{il}(x)l(Y^\varepsilon)^{-1}Y_k^\varepsilon Y_l^\varepsilon\, ds\right]$$

$$= f(x) + \frac{t}{2}\int_{|y|\leq\varepsilon} \sum_{i,j,k,l} \partial_{ij}f(x)h_{ik}h_{il}(x)y_k y_l \nu(dy) + R_\varepsilon(x),$$

where by Lemma 6.9, we have

$$|R_\varepsilon(x)| \leq C\|f\|_{C_p^4}(1 + |x|^{p+4})t^2$$
$$\times \int_{|y|\leq\varepsilon} |y|^4 l(y)^{-1}(1 + (|y|^2 l(y))^{-1}\lambda_\varepsilon t)^{(p-2)/2})\nu(dy).$$

Furthermore

$$L_{d+1}^{2,\varepsilon}f(x) - \frac{1}{2}\int_{|y|\leq\varepsilon}\sum_{i,j,k,l}\partial_{ij}f(x)h_{ik}h_{il}(x)y_k y_l \nu(dy)$$
$$= \sum_{i,j,k,l}\int_{|y|\leq\varepsilon}\int_0^1 (\partial_{ij}f(x+\alpha h(x)y) - \partial_{ij}f(x))\alpha\, d\alpha\, h_{ik}h_{il}(x)y_k y_l \nu(dy).$$

Therefore

$$\left|L_{d+1}^{2,\varepsilon}f(x) - \frac{1}{2}\int_{|y|\leq\varepsilon}\sum_{i,j,k,l}\partial_{ij}f(x)h_{ik}h_{il}(x)y_k y_l \nu(dy)\right|$$
$$\leq C\|f\|_{C_p^4}(1+|x|^{p+3})\int_{|y|\leq\varepsilon}|y|^3 \nu(dy).$$

This finishes the proof. $\square$

In the particular case that $l(y) = y^r$, $r = 2$, the above scheme corresponds to a Asmussen–Rosiński type approach.

The approximation for $L_{d+1}^3$ is defined as follows. Let

$$G_{\varepsilon,l}(dy) = C_{\varepsilon,l}^{-1}l(y)1_{|y|>\varepsilon}\nu(dy),$$
$$C_{\varepsilon,l} = \int_{|y|>\varepsilon}l(y)\nu(dy)$$

and let $Z^{\varepsilon,l} \sim G_{\varepsilon,l}$ and let $S^{\varepsilon,l}$ be a Bernoulli random variable independent of $Z^{\varepsilon,l}$. Then consider the following two cases. If $S^{\varepsilon,l} = 0$ define $\bar{X}_t^{3,\varepsilon}(x) = x$, otherwise $\bar{X}_t^{3,\varepsilon}(x) = x + h(x)l(Z^{\varepsilon,l})^{-1}Z^{\varepsilon,l}$. Then we have the following results.

LEMMA 6.11. *Assume that for $p \geq 2$, $\sup_{\varepsilon \in (0,1]}\int_{|y|>\varepsilon}l(y)^{-p}|y|^{p+1}\nu(dy) < \infty$ and $C_{\varepsilon,l}^{-1}P[S^\varepsilon = 1] \leq Ct$ then Assumption ($\mathcal{M}$) is satisfied with*

$$E[|\bar{X}_{d+1}^{3,\varepsilon}(x)|^p] \leq (1+Kt)|x|^p + K't.$$

PROOF. The result follows clearly from ($f(x) = |x|^p$)

$$P[S^\varepsilon = 1]|E[f(x + h(x)l(Z^{\varepsilon,l})^{-1}Z^{\varepsilon,l}) - f(x)]|$$
$$= C_{\varepsilon,l}^{-1}P[S^\varepsilon = 1]\int_{|y|>\varepsilon}(f(x+h(x)l(y)^{-1}y) - f(x))l(y)\nu(dy)$$
$$\leq Ct(1+|x|^p)\left(1 + \int_{|y|>\varepsilon}l(y)^{-p}|y|^{p+1}\nu(dy)\right). \quad \square$$

LEMMA 6.12. *Assume that for* $f \in C_p^2$, *we have that* $\int_{|y|>\varepsilon} |y|^2 (l(y)^{-1} - 1) + |y|^{p+2}|l(y)^{-1} - 1|^{p+1}\nu(dy) \leq C$ *and* $|C_{\varepsilon,l}^{-1} P[S^{\varepsilon,l} = 1] - t| \leq Ct^2$ *then*

$$|E[f(\bar{X}_t^{3,\varepsilon})] - f(x) - tL_{d+1}^3 f(x)| \leq Ct^2 \|f\|_{C_p^2}(1 + |x|^{p+2}).$$

PROOF. As before let $f \in C_p^2$ then

$$E[f(\bar{X}_t^{3,\varepsilon})] = f(x) + E[f(x + h(x)l(Z^{\varepsilon,l})^{-1}Z^{\varepsilon,l}) - f(x); S^{\varepsilon,l} = 1]$$

$$= f(x) + \int_{|y|>\varepsilon} (f(x + h(x)l(y)^{-1}y) - f(x))l(y)\nu(dy)$$

$$\times C_{\varepsilon,l}^{-1} P[S^{\varepsilon,l} = 1].$$

Then we clearly have that

$$|E[f(\bar{X}_t^{3,\varepsilon})] - f(x) - tL_{d+1}^3 f(x)|$$

$$\leq t \left| \int_{|y|>\varepsilon} \int_0^1 \sum_i (\partial_i f(x + \alpha h(x)l(y)^{-1}y) - \partial_i f(x + \alpha h(x)y))\, d\alpha\, h(x) y \nu(dy) \right|$$

$$\times C_{\varepsilon,l}^{-1} P[S^{\varepsilon,l} = 1] + \left| \int_{|y|>\varepsilon} f(x + h(x)y) - f(x)\nu(dy) \right|$$

$$\times |C_{\varepsilon,l}^{-1} P[S^{\varepsilon,l} = 1] - t|$$

$$\leq C\|f\|_{C_p^2}(1 + |x|^{p+2})t^2.$$

This finishes the proof. □

Using the previous results we can propose various schemes of approximation of order 1 as in Theorem 6.5. We state the simplest type of approximation.

THEOREM 6.6. *Assume that $V_0$, $V$ and $h$ are infinitely differentiable functions with bounded derivatives with $\int_{\mathbf{R}_0^d}(1 \wedge |y|^p)\nu(dy) < \infty$ for all $p \in \mathbf{N}$. Define $\varepsilon \equiv \varepsilon(T,n)$ so that the conditions on Lemmas 6.9, 6.10, 6.11 and 6.12 are satisfied for $t = T/n$ and for appropriate localization functions. Let $\bar{X}_{i,t}^j(x), i = 0, \ldots, d,\ j = 1, \ldots, n$, $n$ copies of the Euler–Maruyama method for $X_{i,t}(x)$.*

*Also, let $\bar{X}_{d+1,T/n}^{i,\varepsilon,j}$, $i = 1, 2, 3$, $j = 1, \ldots, n$, be $n$ independent copies of the schemes defined above. Then the following scheme, $X_T^{(n)} = Y_n^n \circ Y_n^{n-1} \circ \cdots \circ Y_n^1(x)$, $Y_n^j(x) = \bar{X}_{0,T/n}^j \circ \cdots \circ \bar{X}_{d,T/n}^j \circ \bar{X}_{d+1,T/n}^{1,\varepsilon,j} \circ \bar{X}_{d+1,T/n}^{2,\varepsilon,j} \circ \bar{X}_{d+1,T/n}^{3,\varepsilon,j}(x)$ is a first-order discretization scheme.*

Achieving higher-order schemes for the approximation of $L_{d+1}^2$ can be easily obtained from the proof of Lemma 6.10. In fact, the required conditions are as follows. Assume that for $p \geq 2$,

$$\text{(6.15)} \quad \int_{|y|\leq \varepsilon} |y|^4 l(y)^{-1}(1 + (|y|^2 l(y)^{-1}\lambda_\varepsilon t)^{(p-2)/2})\nu(dy) \leq Ct,$$

$$\text{(6.16)} \quad \int_{|y|\leq \varepsilon} |y|^3 \nu(dy) \leq Ct^2.$$

For $L_{d+1}^3$, the idea used in the previous scheme is that the probability of having more than one jump in an interval of size $T/n$ is of order $(T/n)^2$ and therefore they can be neglected if the goal is to achieve a scheme of order 1. Obviously, in order to obtain a higher-order scheme, one has to consider the possibility of more jumps per interval. As an example, we consider the case of at most two jumps per interval with localization $l \equiv 1$.

For $L_{d+1}^3$ one can do the following: let $G_\varepsilon(dy) = C_\varepsilon^{-1} 1_{|y|>\varepsilon}\nu(dy)$, $C_\varepsilon = \int_{|y|>\varepsilon} \nu(dy)$ and let $Z_1^\varepsilon, Z_2^\varepsilon \sim G_\varepsilon$ independent between themselves and let $S_1^\varepsilon$ and $S_2^\varepsilon$ be two independent Bernoulli random variable independent of $Z_1^\varepsilon$, $Z_2^\varepsilon$. Then consider the following cases. If $S_1^\varepsilon = 0$ define $\bar{X}_t^{3,\varepsilon}(x) = x$, if $S_1^\varepsilon = 1$ and $S_2^\varepsilon = 0$ then $\bar{X}_t^{3,\varepsilon}(x) = x + h(x)Z_1^\varepsilon$ and finally if $S_1^\varepsilon = 1$ and $S_2^\varepsilon = 1$ then $\hat{X}_t^{3,\varepsilon}(x) = x + h(x)Z_1^\varepsilon + h(x + h(x)Z_1^\varepsilon)Z_2^\varepsilon$.

Define

$$p_\varepsilon = P[S_1^\varepsilon = 1](1 + P[S_2^\varepsilon = 1]),$$
$$q_\varepsilon = P[S_1^\varepsilon = 1]P[S_2^\varepsilon = 1].$$

In this case we have the following lemma.

LEMMA 6.13. *If $C_\varepsilon^{-1} P[S_1^\varepsilon = 1, S_2^\varepsilon = 0] \leq Ct$ and $C_\varepsilon^{-2} P[S_1^\varepsilon = 1, S_2^\varepsilon = 1] \leq Ct$ then Assumption (M) is satisfied with*

$$E[|\hat{X}_{d+1}^{3,\varepsilon}(x)|^p] \leq (1 + Kt)|x|^p + K't$$

*for all $p \geq 2$.*

PROOF. The result follows clearly from ($f(x) = |x|^p$)

$$P[S_1^\varepsilon = 1, S_2^\varepsilon = 0]|E[f(x + h(x)Z^\varepsilon) - f(x)]|$$
$$\leq Ct(1 + |x|^p)\left(1 + \int_{|y|>\varepsilon} |y|^p \nu(dy)\right),$$
$$P[S_1^\varepsilon = 1, S_2^\varepsilon = 1]|E[f(x + h(x)Z_1^\varepsilon + h(x + h(x)Z_1^\varepsilon)Z_2^\varepsilon) - f(x)]|$$
$$\leq Ct(1 + |x|^p)\left(1 + \left(\int_{|y|>\varepsilon} |y|^p \nu(dy)\right)^2\right). \quad \square$$

LEMMA 6.14.  *Assume that* $|C_\varepsilon^{-1} p_\varepsilon - t| \leq Ct^3$ *and* $|2C_\varepsilon^{-2} q_\varepsilon - t^2| \leq Ct^3$ *then*

$$\left| E[f(\hat{X}_t^{3,\varepsilon})] - f(x) - tL_{d+1}^3 f(x) - \frac{t^2}{2}(L_{d+1}^3)^2 f(x) \right|$$

$$\leq Ct^3 \|f\|_{C_p^2}(1 + |x|^{p+2})\left(1 + \left(\int_{|y|>\varepsilon} |y|\nu(dy)\right)^2\right).$$

PROOF.  As before let $f \in C_p^2$ then

$E[f(\hat{X}_t^{3,\varepsilon})]$

$= f(x) + \int_{|y|>\varepsilon} (f(x + h(x)y) - f(x))\nu(dy) C_\varepsilon^{-1} P[S_1^\varepsilon = 1, S_2^\varepsilon = 0]$

$\quad + E\left[\int_{|y|>\varepsilon} f(x + h(x)y + h(x + h(x)y)Z_2^\varepsilon) - f(x)\nu(dy)\right]$

$\quad \times C_\varepsilon^{-1} P[S_1^\varepsilon = 1, S_2^\varepsilon = 1]$

$= f(x) + L_{d+1}^3 f(x) C_\varepsilon^{-1} P[S_1^\varepsilon = 1, S_2^\varepsilon = 0]$

$\quad + \int_{|y|>\varepsilon}\int_{|y|>\varepsilon} f(x + h(x)y + h(x + h(x)y)y_1) - f(x)\nu(dy)\nu(dy_1)$

$\quad \times C_\varepsilon^{-2} P[S_1^\varepsilon = 1, S_2^\varepsilon = 1]$

$= f(x) + L_{d+1}^3 f(x) C_\varepsilon^{-1}(P[S_1^\varepsilon = 1] + P[S_1^\varepsilon = 1, S_2^\varepsilon = 1])$

$\quad + (L_{d+1}^3)^2 f(x) C_\varepsilon^{-2} P[S_1^\varepsilon = 1, S_2^\varepsilon = 1].$

Therefore

$$\left| E[f(\hat{X}_t^{3,\varepsilon})] - f(x) - tL_{d+1}^3 f(x) - \frac{t^2}{2}(L_{d+1}^3)^2 f(x) \right|$$

$$\leq |L_{d+1}^3 f(x)||C_\varepsilon^{-1} p_\varepsilon - t| + |(L_{d+1}^3)^2 f(x)|\left|C_\varepsilon^{-2} q_\varepsilon - \frac{t^2}{2}\right|.$$

Finally note that

$(L_{d+1}^3)^2 f(x)$

$= \int_{\varepsilon<|y|} L_{d+1}^3 f(x + h(x)y) - L_{d+1}^3 f(x)\nu(dy)$

$= \int_{\varepsilon<|y|}\int_{\varepsilon<|y_1|} (f(x + h(x)y + h(x + h(x)y)y_1)$

$\qquad\qquad - 2f(x + h(x)y) + f(x))\nu(dy_1)\nu(dy)$

$$= \int_{\varepsilon<|y|}\int_{\varepsilon<|y_1|}\int_0^1 \nabla f(x+h(x)y+\alpha h(x+h(x)y)y_1)h(x+h(x)y)y_1$$
$$- \nabla f(x+\alpha h(x)y)h(x)y\,d\alpha\,\nu(dy_1)\nu(dy)$$
$$= \int_{\varepsilon<|y|}\int_{\varepsilon<|y_1|}\int_0^1 \nabla f(x+h(x)y_1+\alpha h(x+h(x)y_1)y)$$
$$\times \int_0^1 \nabla h(x+\beta h(x)y_1)h(x)y_1\,d\beta\,y\,d\alpha\,\nu(dy_1)\nu(dy)$$
$$+ \int_{\varepsilon<|y|}\int_{\varepsilon<|y_1|}\int_0^1\int_0^1 D^2 f(x+\alpha h(x)y$$
$$+\beta(h(x)y_1+\alpha(h(x+h(x)y_1)-h(x))y))$$
$$\times \Big[h(x)y_1$$
$$+\alpha\Big(\int_0^1 \nabla h(x+\gamma h(x)y_1)\,d\gamma\,h(x)y_1\Big)y, h(x)y\Big]\,d\beta\,d\alpha\,\nu(dy_1)\nu(dy).$$

This finishes the proof. □

A similar statement can be achieved if we limit the number of jumps in any interval. The parallel of Theorem 6.6 can also be stated in this case.

6.5. *Example: tempered stable Lévy measure.* Now we consider the previous approximation in the case that the Lévy measure $\nu$ defined on $\mathbf{R}_0$ is given by

$$\nu(dy) = \frac{1}{|y|^{1+\alpha}}(c_+ e^{-\lambda_+|y|}1_{y>0}+c_- e^{-\lambda_-|y|}1_{y<0})\,dy.$$

The Lévy process associated with no Brownian term and the above Lévy measure $\nu$ is called by:

- Gamma: $\lambda_+, c_+ > 0$, $c_- = 0$, $\alpha = 0$.
- Variance gamma: $\lambda_+, \lambda_-, c_+, c_- > 0$, $\alpha = 0$.
- Tempered stable: $\lambda_+, \lambda_-, c_+, c_- > 0$, $0 < \alpha < 2$.

Then, we have that for $\alpha \in [0,1)$

$$\int_{|y|\leq \varepsilon}|y|^k \nu(dy) \sim \varepsilon^{k-\alpha}, \qquad k \geq 1.$$

Then $\sup_{\varepsilon \in (0,1]} \int_{|y|\leq \varepsilon}|y|\nu(dy) < \infty$. For $L^2_{d+1}$, we consider as localization function $l(y) = |y|^r$, then the conditions of Lemma 6.10 are satisfied if $\alpha < r \leq 2$ and $\varepsilon = t^{1/(3-\alpha)}$.

For $L_{d+1}^3$, we consider as localization $l(y) \equiv 1$, then Lemma 6.12 is satisfied, for example, in the following case. Let $P[S^\varepsilon = 1] = e^{-C_\varepsilon a(\varepsilon,t)}$ where $C_\varepsilon \sim \varepsilon^{-\alpha}$, $a(\varepsilon, t) = -\varepsilon^\alpha \log((t^2+t)\varepsilon^{-\alpha})$ as $\varepsilon = t^{1/(3-\alpha)}$ then we have that

$$a = -t^{\alpha/(3-\alpha)} \log((t+1)t^{(3-2\alpha)/(3-\alpha)}).$$

In the case of Lemma 6.14, one choice of parameters is

$$P[S_1^\varepsilon = 1] = t^{(6-3\alpha)/(3-\alpha)}(t+1)(1+t^{\alpha/(3-\alpha)}),$$
$$P[S_2^\varepsilon = 1] = \frac{1}{2(1+t^{\alpha/(3-\alpha)})}.$$

The choice of $r$ in the above scheme is related with variance/importance sampling issues.

*Final comment*: In this article we have presented a general setup to handle what maybe called operator decomposition methods. In particular, the method is useful when considering approximations of expectations of functionals of diffusions (for another similar approach, see Alfonsi [1]). The approximation problem is divided in components, each one driven by a single process. This single process, called the coordinate process can be approximated to a high order using an appropriate (stochastic) Runge–Kutta scheme if the driving process is the Brownian motion. In the case that the driving process is a Lévy process one can decompose the Lévy measure in various pieces to facilitate the analysis. Note that sometimes is not needed to know how to simulate $Y$ but only the functional form of the Lévy measure. In comparison with the proposal presented in [9], where high-order multiple integrals driven by different Wiener processes have to be simulated at each step, we believe that the present methodology is a better scheme.

The issue that local approximations of high order are interesting to study in comparison with Romberg extrapolations as introduced in [19] is similar to the discussion of using Runge–Kutta approximations in comparison with Romberg extrapolations to approximate solutions of ordinary differential equations. We believe this article helps to open the path in this direction. In fact, it is somewhat clear from Theorem 4.2 that the leading constants in a Euler+Romberg method and a Runge–Kutta method do not coincide.

Finally, we used the structure of this construction to easily introduce and analyze the asymptotic error of an approximating scheme for solutions of stochastic differential equations driven by Lévy processes with possibly infinite activity.

## APPENDIX

In this section we assume condition $(\mathcal{M}_P)$.

LEMMA A.1. *Let $f = f_s(x) \in C_p^{1,2}([0,T] \times \mathbf{R}^N)$. Then a map $s \mapsto P_s f_s(x)$ is Lipschitz continuous for all $x \in \mathbf{R}^N$.*

PROOF. Note that

$$|P_t f_t(x) - P_s f_s(x)| \leq |P_t f_t(x) - P_t f_s(x)| + |P_t f_s(x) - P_s f_s(x)|.$$

Using the Lipschitz properties of $t \mapsto f_t(x)$ and $t \mapsto P_t f_s(x)$, the proof follows. □

LEMMA A.2. *Let $g:[0,T] \to \mathbf{R}$ be a Lipschitz continuous function. Then we have*

(A.1) $$\left| \frac{T}{n} \sum_{k=1}^{n} g(kT/n) - \int_0^T g(s)\,ds \right| \leq \frac{C(T,g)}{n}.$$

PROOF. From the assumption we immediately obtain

$$\left| \frac{T}{n} g(kT/n) - \int_{(k-1)T/n}^{kT/n} g(s)\,ds \right| \leq \frac{C}{n^2},$$

where $C$ depends on $T$ and the Lipschitz coefficient of $g$. This implies (A.1). □

**A.1. Some properties of Lévy driven SDEs.** We start with the differentiability properties of $X_t(x)$ in $x$. The following material can be found in [7, 8, 10, 17] and [18]. We quote them here for completeness.

LEMMA A.3. *There exists a version of $X_t(x)$ such that a map $x \mapsto X_t(x)$ is infinite times continuously differentiable almost surely and in the $L^p$-sense. Moreover, we have for $p \geq 2$,*

(A.2) $$E\left[ \sup_{0 \leq t \leq T} |X_t(x)|^p \right] \leq C(p,T)(1 + |x|^p)$$

*and*

(A.3) $$\sup_{x \in \mathbf{R}^N} E\left[ \sup_{0 \leq t \leq T} |\partial_x^\alpha X_t(x)|^p \right] < \infty$$

*for any multi-index $\alpha$ with $|\alpha| \geq 1$.*

PROPOSITION A.4. *Let $f \in C_p^m$ with $p \geq 2$.*

(i) *Then $P_t f \in C_p^m$ for all $t \geq 0$ and*

(A.4) $$\sup_{t \in [0,T]} \|P_t f\|_{C_p^m} \leq C \|f\|_{C_p^m}.$$

(ii) If $m \geq 2$, then $Lf \in C_{p+2}^{m-2}$ and
$$\|Lf\|_{C_{p+2}^{m-2}} \leq C\|f\|_{C_p^m}.$$

(iii) If $f \in C_p^{1,m}([0,T] \times \mathbf{R}^N)$, then $(\partial_t Lf)(t,x) = (L\partial_t f)(t,x)$.

PROOF. The proof of (i) follows by interchange of derivation and expectation together with the moment estimates in Lemma A.3. Recall that $L = \sum_{i=0}^{d+1} L_i$ as defined in (2.5).

(ii) We only do the proof for $L_{d+1}$. We have
$$\left|\int (f(x+h(x)y) - f(x) - \nabla f(x)h(x)\tau(y))\nu(dy)\right|$$
$$\leq \left|\int \nabla f(x)h(x)(y - \tau(y))\nu(dy)\right| + \left|\int\int_0^1 \frac{d^2}{d\theta^2} f(x + \theta h(x)y)\,d\theta\,\nu(dy)\right|$$
$$\leq C\|f\|_{C_p^m}(1 + |x|^{p+2}). \qquad \square$$

PROPOSITION A.5. *Let $f \in C_p^2$. Then $P_t$ and $L$ are commutative and $u_f(t,x) := P_t f(x)$ is the solution of the integro-differential equation:*
$$\begin{cases} \dfrac{d}{dt} u_f(t,x) = L u_f(t,x), \\ u_f(0,x) = f(x). \end{cases}$$

Let $f \in C_p^{2m+2}$. Then the commutativity of $P_t$ and $L$ implies that $L^m u_f$ ($= u_{L^m f}$) is differentiable in $t$ and the solution to similar integro-differential equations. That is,
$$\begin{cases} \dfrac{d}{dt}(L^m u_f)(t,x) = L(L^m u_f)(t,x), \\ (L^m u_f)(0,x) = (L^m f)(x) \end{cases}$$
for each $m \geq 0$. Consequently, applying Taylor's expansion to $u_f$, we have:

PROPOSITION A.6. *For $f \in C_p^{2m+2}$,*
$$P_t f(x) = \sum_{k=0}^m \frac{t^k}{k!} L^k f(x) + \int_0^t \frac{(t-s)^m}{m!} P_s(L^{m+1} f)(x)\,ds.$$

*Furthermore, if $f \in C_p^m$ with $m \geq 2$. Then $P_t f \in C_{p+2}^{1,m-2}$.*

Summarizing this section, we have

COROLLARY A.7.   $P_t f(x) = E[f(X_t(x))]$ and $Q_t^i f(x) = E[f(X_t^i(x))]$ ($i = 0, 1, \ldots, d+1$) *satisfy conditions* $(\mathcal{M})$ *and* $\mathcal{R}(m, t^m)$. *That is, for* $p \in \mathbf{N}$,

$$E[|X_t(x)|^{2p}] \leq (1 + Kt)|x|^{2p} + K't$$

*for some constant* $K = K(T, p)$, $K' = K'(T, p) > 0$ *and*

$$J_{\leq m}(P_t) = \sum_{k=0}^{m} \frac{t^k}{k!} L^k,$$

$$J_{\leq m}(Q_t^i) = \sum_{k=0}^{m} \frac{t^k}{k!} L_i^k$$

*for any* $m \in \mathbf{N}$.

MITSUBISHI UFJ TRUST INVESTMENT  
TECHNOLOGY INSTITUTE CO., LTD.  
5-6, SHIBA 2-CHOME  
MINATO-KU  
TOKYO 105-0014  
JAPAN  
E-MAIL: tanaka@mtec-institute.co.jp

OSAKA UNIVERSITY  
GRADUATE SCHOOL OF ENGINEERING SCIENCES  
MACHIKANEYAMA CHO 1-3  
OSAKA 560-8531  
JAPAN  
E-MAIL: kohatsu@sigmath.es.osaka-u.ac.jp